\documentclass{amsart}

\usepackage[utf8]{inputenc}
\usepackage[english]{babel}
\usepackage[T1]{fontenc}

\usepackage{arydshln}
\usepackage{amsmath}
\usepackage{amssymb}
\usepackage{amsthm}
\usepackage{mathabx}
\usepackage{tikz}
\usepackage{tikz-cd}
\usepackage{todonotes}
\usepackage{hyperref}
\usepackage{cleveref}
\usepackage{footmisc}
\usepackage{bbm}
\usepackage{afterpage}
\usepackage{mathrsfs}
\usepackage{sseq}
\usepackage{spectralsequences}
\usepackage{extarrows}
\usepackage{enumitem}
\usepackage{stmaryrd}
\usepackage{pdflscape}
\hypersetup{
  colorlinks   = true, 
  urlcolor     = blue, 
  linkcolor    = blue, 
  citecolor   = blue 
}

\newcommand{\bp}{\mathrm{bp}}

\DeclareMathOperator{\ANSS}{ANSS}

\newcommand{\BP}{\mathrm{BP}}

\DeclareMathOperator{\KO}{KO}
\DeclareMathOperator{\ko}{ko}

\DeclareMathOperator{\MU}{MU}

\newcommand{\Sph}{\mathbf{S}}

\newcommand{\TMF}{\mathsf{TMF}}
\newcommand{\tmf}{\mathsf{tmf}}
\newcommand{\Tmf}{\mathsf{Tmf}}

\DeclareMathOperator{\Fun}{Fun}
\DeclareMathOperator{\Fil}{Fil}
\DeclareMathOperator{\map}{map}

\DeclareMathOperator{\Sp}{Sp}

\DeclareMathOperator{\Syn}{Syn}

\DeclareMathOperator{\fib}{fib}
\DeclareMathOperator{\coker}{coker}
\newcommand{\id}{\mathrm{id}}

\newcommand{\op}{\mathrm{op}}

\newcommand{\E}{\mathbf{E}}

\newcommand{\F}{\mathbf{F}}

\newcommand{\J}{\mathsf{J}}
\renewcommand{\j}{\mathrm{j}}

\newcommand{\Q}{\mathsf{Q}}

\newcommand{\gX}{\mathbf{X}}
\newcommand{\gY}{\mathbf{Y}}

\newcommand{\Z}{\mathbf{Z}}

\newcommand{\al}{\alpha}
\newcommand{\be}{\beta}

\newcommand{\Ga}{\Gamma}
\newcommand{\eps}{\varepsilon}
\newcommand{\kappabar}{\bar{\kappa}}

\let\oldbold\boldsymbol
\renewcommand{\boldsymbol}[1]{\textcolor{blue}{\oldbold{#1}}}

\usepackage{mathtools}

\theoremstyle{theorem}\numberwithin{equation}{section}
\newtheorem{theorem}[equation]{Theorem}
\crefname{theorem}{{th}.\negthickspace}{{ths}.\negthickspace}
\Crefname{theorem}{{Th}.\negthickspace}{{Ths}.\negthickspace}
\newtheorem{theoremalph}{Theorem}

\newtheorem{coralph}[theoremalph]{Corollary}

\crefname{theoremalph}{{th}.\negthickspace}{{ths}.\negthickspace}
\Crefname{theoremalph}{{Th}.\negthickspace}{{Ths}.\negthickspace}

\Crefname{problem}{{Prb}.\negthickspace}{{Prbs}.\negthickspace}
\newtheorem{prop}[equation]{Proposition}
\Crefname{prop}{{Pr}.\negthickspace}{{Prs}.\negthickspace}
\newtheorem{lemma}[equation]{Lemma}
\Crefname{lemma}{{Lm}.\negthickspace}{{Lms}.\negthickspace}
\newtheorem{cor}[equation]{Corollary}
\Crefname{cor}{{Cor}.\negthickspace}{{Cors}.\negthickspace}

\Crefname{conjecture}{{Conj}.\negthickspace}{{Conjs}.\negthickspace}

\theoremstyle{definition}\numberwithin{equation}{section}
\newtheorem{mydef}[equation]{Definition}
\Crefname{mydef}{{Df}.\negthickspace}{{Dfs}.\negthickspace}

\Crefname{variant}{{Var}.\negthickspace}{{Vars}.\negthickspace}

\Crefname{recall}{{Rcl}.\negthickspace}{{Rcls}.\negthickspace}

\Crefname{construction}{{Con}.\negthickspace}{{Cons}.\negthickspace}

\Crefname{ass}{{As}.\negthickspace}{{As}.\negthickspace}
\newtheorem{question}[equation]{Question}
\Crefname{question}{{Q}.\negthickspace}{{Qs}.\negthickspace}
\newtheorem{notation}[equation]{Notation}
\Crefname{notation}{{Nt}.\negthickspace}{{Nts}.\negthickspace}

\Crefname{situation}{{St}.\negthickspace}{{Sts}.\negthickspace}

\theoremstyle{remark}\numberwithin{equation}{section}
\newtheorem{example}[equation]{Example}
\Crefname{example}{{Ex}.\negthickspace}{{Exs}.\negthickspace}

\Crefname{nonexample}{{NonEx}.\negthickspace}{{NonEx}.\negthickspace}

\Crefname{claim}{{Clm}.\negthickspace}{{Clms}.\negthickspace}
\newtheorem{remark}[equation]{Remark}
\Crefname{remark}{{Rmk}.\negthickspace}{{Rmks}.\negthickspace}

\Crefname{idea}{{Id}.\negthickspace}{{Ids}.\negthickspace}
\newtheorem{warn}[equation]{Warning}
\Crefname{warn}{{Warn}.\negthickspace}{{Warns}.\negthickspace}

\Crefname{figure}{{Fig.}\negthickspace}{{Figs.}\negthickspace}
\Crefname{footnote}{{Fn.}\negthickspace}{{Fn.}\negthickspace}

\Crefname{part}{{\textsection}\negthickspace}{{\textsection}\negthickspace}
\Crefname{chapter}{{\textsection}\negthickspace}{{\textsection}\negthickspace}
\Crefname{section}{{\textsection}\negthickspace}{{\textsection}\negthickspace}
\Crefname{subsection}{{\textsection}\negthickspace}{{\textsection}\negthickspace}
\Crefname{appendix}{{\textsection}\negthickspace}{{\textsection}\negthickspace}

\DeclarePairedDelimiter{\angbr}{\langle}{\rangle}

\begin{sseqdata}[name = koexample1, classes = fill, Adams grading, tick step = 2,
major tick step = 2,
minor tick step = 2]

\class[rectangle](0,0)\class[rectangle](4,0)
\class(0,4)\class(-1,7)
\class(1,1)\class(2,2)\class(3,3)\class(4,4)\class(5,5)\class(6,6)\class(7,7)
\class(5,1)\class(6,2)\class(7,3)\class(8,4)
\structline(0,0)(7,7)\structline(4,0)(8,4)
\classoptions["{u^2}"  {left, pin = black}](4,0)
\classoptions["{\eta}"  {above left, pin = black}](1,1)
\end{sseqdata}

\begin{sseqdata}[name = koexample2, Adams grading, tick step = 2,
major tick step = 2,
minor tick step = 2]

\class[rectangle, blue](0,0)
\class[blue](0,4)\class[blue](-1,7)
\class[blue](1,1)\class[blue](2,2)\class[blue](3,3)\class[blue](4,4)\class[blue](5,5)\class[blue](6,6)\class[blue](7,7)
\class[blue](5,1)\class[blue](6,2)\class[blue](7,3)\class[blue](8,4)
\structline(0,0)(7,7)\structline(5,1)(8,4)
\classoptions["{\eta}" {below right, pin = black}](1,1)

\class[red](0,2)\class[red](1,3)\class[red](2,4)\class[red](3,5)\class[red](4,6)\class[red](5,7)
\structline(0,2)(5,7)
\class[diamond, red, "8"](3,1)\class[red](4,2)\class[red](5,3)\class[red](6,4)
\structline(3,1)(6,4)
\end{sseqdata}

\begin{document}
\title{
On periodic families in the stable stems of height two
}

\author{Christian Carrick}
\email{carrick@math.uni-bonn.de}
\author{Jack Morgan Davies}
\email{davies@math.uni-bonn.de}

\date{\today}

\maketitle
\begin{abstract}
    We discover a host of infinite periodic families in the $2$-primary stable homotopy groups of spheres. We also confirm the existence of many families predicted by Hopkins--Mahowald. These families appear in nineteen different congruence classes of degrees modulo $192$, seven of them consist of simple $4$-torsion elements, and another four of simple $8$-torsion. They all vanish in the homotopy groups of the spectrum $\TMF$ of topological modular forms, but we show that they are detected in the fixed-points of $\TMF$ with respect to an Atkin--Lehner involution. As a consequence, we confirm the existence of exotic spheres in all dimensions congruent to $72$, $144$, and $168$ modulo $192$.
\end{abstract}    

\setcounter{tocdepth}{1}
\tableofcontents
\vspace{-1.5cm}
\section{Introduction}
The basic building blocks of the stable homotopy category are the stable homotopy groups of spheres $\pi_\ast\Sph$. Due to the abundance of rich arithmetic and geometric information contained within these groups, explicit computations are both difficult and highly sought after. Modern stem-by-stem computations have taken large strides recently, see \cite{iwxpisph} for a survey and \cite{linwangxu} for a culmination of these efforts. Despite these monumental advances, it remains incredibly difficult to answer basic questions, such as if a certain family of elements defined in the Adams spectral sequence (ASS), one of our most powerful tools to compute these groups, survives to represent classes in $\pi_\ast\Sph$; see Adams \cite{adamshopf} and Hill--Hopkins--Ravenel \cite{hhr} for examples of such results and their wider applications.

\emph{Chromatic homotopy theory} is a widely successful program to organise $\pi_\ast\Sph$ into \emph{periodic families}. The first glimpse of this behaviour comes from Adams \cite{adamsjofx}, who shows that $\pi_\ast\Sph$ contains six $8$-periodic families of elements detected by a combination of the $8$-fold Bott periodicity of \emph{real topological $K$-theory} $\KO$, the \emph{Adams operations} $\psi^k$, and the \emph{$J$-homomorphism}. From the perspective of chromatic homotopy theory, these families are \emph{$v_1$-periodic} families, meaning that they are detected by classes in the $1$-line of the Adams--Novikov spectral sequence (ANSS), a variant of the ASS due to Novikov \cite{anoneline}. The first $v_2$-periodic families were later discovered by Smith \cite{smithbeta} and Miller--Ravenel--Wilson \cite{MRW77},\footnote{In \cite{MRW77}, Miller--Ravenel--Wilson also discover the fundamental $v_3$-periodic family at primes $p\geq 7$. Virtually nothing else is known about $v_3$-periodic families nor about $v_h$-periodic families for $h\geq 4$.} although these results only produce such families in the $p$-completion $\pi_\ast\Sph_p$ at a prime number $p\geq 5$. As demonstrated by \cite{adamshopf,hhr}, much of the geometric information in $\pi_\ast\Sph$ seems to be concentrated at the prime $2$. 

Historically, much of the difficulty of the primes $2$ and $3$ stemmed from the lack of a $v_2$-periodic analogue of $\KO$. In seminal work, Hopkins--Mahowald \cite{hopkinsfirsttmficm} discovered this $v_2$-periodic analogue, originally called a higher real $K$-theory, and now known as the spectrum $\TMF$ of \emph{topological modular forms}. Similar to $\KO$, the cohomology theory $\TMF$ has connections to differential geometry and physics, and as the name suggests, also to number theory and arithmetic geometry. Hopkins--Mahowald then used $\TMF$ to identify a vast array of $v_2$-periodic families at the primes $2$ and $3$, results which were later confirmed by Behrens--Mahowald--Quigley \cite{hurewicztmf} and Belmont--Shimomura \cite{tmfthree}. 

All of the $v_2$-periodic families of \cite{hurewicztmf,tmfthree} are \emph{detected} by $\TMF$, meaning that they are nonzero in the image of the unit map $\pi_\ast\Sph \to \pi_\ast\TMF$, also known as the \emph{Hurewicz image} of $\TMF$. Some of Adams' $v_1$-periodic families have nonzero image in $\KO$, such as those generated by the first Hopf classes $\eta$ and $\eta^2$, but some do not, such as the family generated by $\eta^3$.
Producing $v_2$-periodic families which are not detected by $\TMF$ is difficult, partly due to the lack of a $v_2$-analogue of the $J$-homomorphism, which Carmeli \cite{carmeli_nohigherJ} recently showed cannot exist. Until now, the only known $v_2$-periodic families at the prime $2$ which vanish in $\pi_\ast\TMF$ are the two families discovered by Behrens--Hill--Hopkins--Mahowald \cite{exoticwithcokerJ} and the nine families discovered by Bhattacharya--Bobkova--Quigley \cite{bbq,bq}, all with the help of $\TMF$.\footnote{In \cite{heighttwojat3}, we discovered various $v_2$-periodic families at the prime $3$ not detected by $\TMF$.}

\subsection{Main result}
In this article, we use operations on $\TMF$ to produce a host of new, nontrivial $v_2$-periodic families in $\pi_*\Sph$ that vanish in $\TMF$. Our methods also systematically catalogue all such previously known $v_2$-periodic families.

\begin{theoremalph}\label{mainthm:nonzerofamilies}
    \Cref{tab:classes} represents $125$ nonzero $v_2^{32}$-periodic\footnote{The superscript in $v_2^{32}$ refers to the fact that these families are precisely $|v_2^{32}| = 192$-periodic.} families in $\pi_d\Sph_2$ with $d$ congruent to the value in the Degree-column modulo $192$, of order indicated by the Group-column, and all with trivial image in $\pi_\ast\TMF$. The generators for every family in a particular row all have the same image in the ANSS for $\TMF$ up to higher filtration, as indicated by the $\ANSS(\TMF)$-column. A full depiction of these generators is displayed in \Cref{list}.
\end{theoremalph}

Our methods apply quite broadly and capture many nonzero $v_2^{32}$-periodic families at once. In particular, we are able to reconfirm the nonvanishing of families mentioned in \cite{tmfbook,exoticwithcokerJ,bbq,bq}, see \Cref{ssec:relatedwork} for more, and we discover 50 additional families not appearing in the literature. From our point of view, these exact numbers are not the main take-away of this article, but rather the efficacy and simplicity of the methods, outlined in \Cref{ssec:proofsection}. The proof of Theorem \ref{mainthm:nonzerofamilies} gives a unified approach to detecting $v_2^{32}$-periodic families that both recovers the known families appearing in the literature and pushes much further to detect a host of new families.

\begin{table}[h]
\centering
{
\begin{tabular}{|c|c|c|c|c|c|c|c|}
\hline
Degree  &Group  &   $\ANSS(\TMF)$                         &   AN-filts.&   \# families     &   $\J_0(3)$       &   Notes               \\  \hline
$23a$   &$\Z/8$ &   $\nu{\kappabar}$                      &   $[5]$    &   $1$             &   $\checkmark$    &   $\dagger$           \\  
$23b$   &$\Z/2$ &   $\eta\eps\kappa$                      &   $[5,7]$  &   $4$             &   $\checkmark$    &   $\dagger$           \\  
$26$    &$\Z/2$ &   ${\nu^2\kappabar}$                    &   $[6]$    &   $2$             &                   &   $\dagger$           \\  \hline

$47a$   &$\Z/4$ &   ${[2\nu\kappabar\Delta]}$             &   $[5]$    &   $1$             &   $\checkmark$    &   $\dagger$           \\  
$47b$   &$\Z/2$ &   $\eta\kappa[\eps\Delta]$              &   $[5,7]$  &   $5$             &   $\checkmark$    &   $\dagger$           \\  
$48$    &$\Z/2$ &   $\nu{[\eta\kappabar\Delta]}$          &   $[6,10]$ &   $5$             &   $\checkmark$    &   $\dagger$           \\  \hline

$71a$   &$\Z/8$ &   ${[\nu\kappabar\Delta^2]}$            &   $[5]$    &   $1$             &   $\checkmark$    &   $\dagger$           \\  
$71b$   &$\Z/2$ &   $[\eps\Delta][\eta\kappa\Delta]$      &   $[5,7]$  &   $7$             &   $\checkmark$    &   $\dagger$           \\  
$72$    &$\Z/2$ &   $\eta[\nu\kappabar\Delta^2]$          &   $[6,10]$ &   $3$             &   $\checkmark$    &   $\dagger$           \\  
$73$    &$\Z/2$ &   $\eta^2[\nu\kappabar\Delta^2]$        &   $[7,11]$ &   $10$            &   $\checkmark$    &   $\dagger$           \\  
$74a$   &$\Z/4$ &   $\nu{[\nu\kappabar\Delta^2]}$         &   $[6]$    &   $3$             &                   &   $\dagger$           \\  
$74b$   &$\Z/2$ &   ${\kappa\kappabar^3}$                 &   $[14]$   &   $2$             &                   &   $\dagger$           \\  \hline

$95$    &$\Z/2$ &  ${[4\nu\kappabar\Delta^3]}$            &   $[5]$    &   $1$             &   $\checkmark$    &   $\dagger$           \\  \hline

$119a$  &$\Z/4$ &   $\nu{[2\kappabar\Delta^4]}$           &   $[5]$    &   $1$             &   $\checkmark$    &                       \\  
$119b$  &$\Z/2$ &   $\eta\eps[\kappa\Delta^4]$            &   $[5,7]$  &   $11$            &   $\checkmark$    &   $\dagger$           \\  
$120$   &$\Z/2$ &   $\nu{[\eta\kappabar\Delta^4]}$        &   $[6,24]$ &   $2$             &   $\checkmark$    &   $\dagger$           \\  
$122$   &$\Z/2$ &   ${\kappabar[\nu^2\Delta^4]}$          &   $[6]$    &   $2$             &                   &                       \\  \hline

$143$   &$\Z/2$ &   $\eps[\eta\kappa\Delta^5]$            &   $[5]$    &   $8$             &   $\checkmark$    &                       \\  
$144$   &$\Z/2$ &   $\kappa\kappabar[\kappa\Delta^4]$     &   $[8,10]$ &   $5$             &   $\checkmark$    &   $\dagger$           \\  
$145$   &$\Z/2$ &   $\nu\kappa[\eps\Delta^5]$             &   $[5,25]$ &   $13$            &   $\checkmark$    &                       \\  \hline

$167$   &$\Z/2$ &   ${[\eta\kappa\Delta][\eps\Delta^5]}$  &   $[5]$    &   $12$            &   $\checkmark$    &   $\dagger$           \\  
$168$   &$\Z/2$ &   ${\kappabar^2[\eps\Delta^5]}$         &   $[10]$   &   $2$             &   $\checkmark$    &                       \\  
$169$   &$\Z/2$ &   $\kappa\kappabar{[\eta\kappa\Delta^5]}$&  $[9,11]$ &   $14$            &   $\checkmark$    &                       \\  
$170a$  &$\Z/8$ &   ${\kappabar[\nu^2\Delta^6]}$          &   $[6]$    &   $2$             &                   &                       \\  
$170b$  &$\Z/4$ &   ${\kappabar^3[\kappa\Delta^4]}$       &   $[14]$   &   $2$             &                   &                       \\  
$170c$  &$\Z/2$ &   $\nu^3[\nu\kappa\Delta^6]$            &   $[6,26]$ &   $6$             &                   &                       \\  \hline
        &       &                                         &            &   $\Sigma = 125$  &                   &                       \\  \hline
\end{tabular}
}
\caption{\small
{
Nonzero $v_2^{32}$-periodic families in $\pi_\ast \Sph_2$ of \Cref{mainthm:nonzerofamilies} in degrees modulo $192$. The Notes-column refers to \Cref{ssec:relatedwork}.
}}\label{tab:classes}
\end{table}

\subsection{Naming the periodic families}\label{ssec:namingperiodicfamilies}
The number of families in \Cref{mainthm:nonzerofamilies} is counted in a \emph{maximal} sense. For example, \Cref{tab:classes} indicates that, aside from the family generated by $4\nu\kappabar$, there are four $2$-torsion infinite $v_2^{32}$-periodic families generated by classes in degree $23$, and that the images of these generators in the ANSS for $\TMF$ are all $\eta\eps\kappa$. In fact, these generators also agree in $\pi_{23}\Sph$, $\eta\eps\kappa = \nu^3\kappa$. To construct these four families, recall that the elements $\nu^2, \eps$, and $\kappa$ are all themselves $v_2^{32}$-periodic, courtesy of \cite{hurewicztmf}, hence we have associated $v_2^{32}$-periodic families denoted as
\begin{equation}\label{eq:basicperiodicfamilies}
    \{\nu^2_n\}, \qquad \{\eps_n\}, \qquad \{\kappa_n\}
\end{equation}
with $x_0 = x$. The four families in row $23b$ are then given by
\begin{equation}\label{eq:compoundfamilies}\{\eta \eps_n \kappa\}, \qquad \{ \eta \eps \kappa_n\}, \qquad \{\nu \nu^2_n \kappa\}, \qquad \{\nu^3 \kappa_n\},\end{equation}
in other words, we \emph{periodify} in each $v_2^{32}$-periodic factor. We count in this way to reflect the four separate facts that the family $\{\eps_n\}$ supports multiplication by $\eta\kappa$, the family $\{\kappa_n\}$ supports multiplication by $\eta\eps$, the family $\{\nu_n^2\}$ supports multiplication by $\nu\kappa$, and the family $\{\kappa_n\}$ supports multiplication by $\nu^3$. These facts do not imply each other because, although the zeroth elements in the families of (\ref{eq:compoundfamilies}) all agree by design, the families need not agree. This is because there are \emph{different} choices hidden in each of the periodic families of (\ref{eq:basicperiodicfamilies}): a choice of type $2$ complex $M$, a $v_2^{32}$-self map of $M$, and a lift of each generator to the top-cell of $M$. A consequence of the \emph{periodicity theorem}, however, is that there exists a positive integer $d$ for each pair of families in (\ref{eq:compoundfamilies}), such that the $n$th entries in this pair of families agree for all $n$ divisible by $d$; see \Cref{cortoinfiniteagreement}.


In \Cref{ssec:countingmethodology}, we methodically justify every family counted in \Cref{tab:classes}, row by row. To summarise the generators, we have \Cref{list}.

\begin{table}[h]
\centering
{
\begin{tabular}{|c|c|}
\hline
Degree  &   Generators                     \\  \hline
$23a$   &   $\nu\boldsymbol{\kappabar}$                                                                 \\  \hdashline
$23b$   &   $\eta\boldsymbol{\eps\kappa},\nu\boldsymbol{\nu^2 \kappa}$                   \\  \hdashline
$26$    &   $\boldsymbol{\nu^2\kappabar}$                    \\  \hline

$47a$   &   $\boldsymbol{[2\nu\kappabar\Delta]}$             \\  \hdashline
$47b$   &   $\eta\boldsymbol{\kappa[\eps\Delta]}, \boldsymbol{\eps[\eta\kappa\Delta]},$\\
&  $\eta^2\boldsymbol{[\eta\kappabar\Delta]} $     \\  \hdashline
$48^\ast$    &   $\nu\boldsymbol{[\eta\kappabar\Delta]}, \boldsymbol{\eps}\kappabar\boldsymbol{\kappabar}$,\\
&  $\eta\boldsymbol{[2\nu\kappabar\Delta]},\kappa\boldsymbol{\kappa\kappabar} $                    \\  \hline

$71a$   &   $\boldsymbol{[\nu\kappabar\Delta^2]}$            \\  \hdashline
$71b$   &   $\boldsymbol{[\eps\Delta][\eta\kappa\Delta]}, \boldsymbol{\nu^2[\nu\kappa\Delta^2]}$,\\
&  $\eta\boldsymbol{[\eta^2\kappabar\Delta^2]}, \nu\boldsymbol{\kappa[\nu^2\Delta^2]}$     \\  \hdashline
$72$    &   $\eta\boldsymbol{[\nu\kappabar\Delta^2]}, \kappabar\boldsymbol{\kappabar[\eps\Delta]}$        \\  \hdashline
$73$    &   $\eta^2\boldsymbol{[\nu\kappabar\Delta^2]},\eta\kappabar\boldsymbol{\kappabar[\eps\Delta]} ,$\\
&  $ \nu\boldsymbol{[\eta^2\kappabar\Delta^2]}, \boldsymbol{\eps\kappabar[\eta\kappabar\Delta]},$\\
&  $  \boldsymbol{\kappa\kappabar[\eta\kappa\Delta]}$\\  \hdashline
$74a$   &   $\nu\boldsymbol{[\nu\kappabar\Delta^2]}$, $\boldsymbol{\kappabar[\nu^2\Delta^2]}$          \\\hdashline
$74b$   &   $\boldsymbol{\kappa}\kappabar^2\boldsymbol{\kappabar}$                 \\  \hline

$95$    &  $\boldsymbol{[4\nu\kappabar\Delta^3]}$          \\  \hline

$119a$  &   $\nu\boldsymbol{[2\kappabar\Delta^4]}$           \\  \hline

$119b$  &   $\eta\boldsymbol{\eps[\kappa\Delta^4]}, \eta\boldsymbol{\kappa[\eps\Delta^4]},$\\  \hline

\end{tabular}
\begin{tabular}{|c|c|}
\hline

(cont.)&  $ \eta^2\boldsymbol{[\eta\kappabar\Delta^4]},\nu\boldsymbol{\nu^2[\kappa\Delta^4]},$\\
&  $ \nu\boldsymbol{\kappa[\nu^2\Delta^4]}, \boldsymbol{[\nu^2\Delta^2][\nu\kappa\Delta^2]}$     \\  \hdashline
$120$   &   $\nu\boldsymbol{[\eta\kappabar\Delta^4]}, \kappabar^5\boldsymbol{\kappabar}$                       \\ \hdashline
$122$   &   $\boldsymbol{\kappabar[\nu^2\Delta^4]}$                       \\  \hline

$143$   &   $\boldsymbol{\eps[\eta\kappa\Delta^5]}, \eta\boldsymbol{\kappa[\eps\Delta^5]},$\\
&  $ \eta\boldsymbol{[\eps\Delta][\kappa\Delta^4]}, \boldsymbol{[\eta\kappa\Delta][\eps\Delta^4]}$            \\  \hdashline
$144$   &   $\boldsymbol{\kappa\kappabar[\kappa\Delta^4]}, \kappabar\boldsymbol{\kappabar[\eps\Delta^4]}$         \\  \hdashline
$145^\ast$   &   $\nu\boldsymbol{\kappa[\eps\Delta^5]},\kappabar^4\boldsymbol{\kappabar[\eta\kappabar\Delta]} ,$\\
&  $ \eta\kappabar\boldsymbol{\kappabar[\eps\Delta^4]},\boldsymbol{\eps\kappabar[\eta\kappabar\Delta^4]},$\\
&  $ \eta\boldsymbol{\kappa\kappabar[\kappa\Delta^4]},\nu\boldsymbol{[\eps\Delta][\kappa\Delta^4]}$  \\  \hline

$167$   &   $\boldsymbol{[\eta\kappa\Delta][\eps\Delta^5]}, \nu\boldsymbol{\kappa[\nu^2\Delta^6]},$\\
&  $ \nu\boldsymbol{[\nu^2\Delta^2][\kappa\Delta^4]},\boldsymbol{[\nu\kappa\Delta^2][\nu^2\Delta^4]},$\\
&  $ \boldsymbol{[\eps\Delta][\eta\kappa\Delta^5]},\boldsymbol{\nu^2[\nu\kappa\Delta^6]} $           \\  \hdashline
$168$   &   $\kappabar\boldsymbol{\kappabar[\eps\Delta^5]}$         \\  \hdashline
$169$   &   $\boldsymbol{\kappa\kappabar[\eta\kappa\Delta^5]},\eta\kappabar\boldsymbol{\kappabar[\eps\Delta^5]} ,$\\
&  $ \boldsymbol{\kappabar[\eta\kappa\Delta][\kappa\Delta^4]},\boldsymbol{\kappabar[\eps\Delta][\eta\kappabar\Delta^4]},$\\
&  $ \boldsymbol{\kappabar[\eta\kappabar\Delta][\eps\Delta^4]}$     \\  \hdashline
$170a$   &   $\boldsymbol{\kappabar[\nu^2\Delta^6]}$\\ \hdashline
$170b$   &   $\kappabar^2\boldsymbol{\kappabar[\kappa\Delta^4]}$\\ \hdashline
$170c$   &   $\nu\boldsymbol{\nu^2[\nu\kappa\Delta^6]},\kappabar^4\boldsymbol{\kappabar[\eta^2\kappabar\Delta^2]} $,\\
&   $\eta^2\kappabar\boldsymbol{\kappabar[\eps\Delta^5]}$\\    \hline
\end{tabular}
}
\caption{\small
{
Summary of $v_2^{32}$-periodic families of \Cref{mainthm:nonzerofamilies} given by periodicity generators.
}}\label{list}
\end{table}

Bold lettering in blue refers to which elements can act as the periodicity generator, so the periodicity generators in row $23b$ are given by $\nu^2$, $\eps$, and $\kappa$, as described above. If multiple factors of a certain periodicity generator occur, we only bold one of them, as the families generated by each of these factors coincide. In general, each bold symbol is known to be $v_2^{32}$-periodic by either \cite{hurewicztmf} or \Cref{lm:47,lm:71,lm:95}. Square brackets indicate indecomposable elements in $\pi_\ast\Sph$ and elements are named by their image in the ANSS for $\TMF$.

Translating \Cref{list} into the $\#$ families-column of \Cref{tab:classes} is done by counting the number of bold factors appearing in each row, with two caveats in degrees $48$ and $145$, where one over-counts by $1$ each time; see \Cref{ssec:countingmethodology} for more details.

\subsection{Consequences for exotic spheres}
Geometric applications of this theorem are immediate. Kervaire--Milnor \cite{kervairemilnor} construct groups $\Theta_k$ of \emph{homotopy $k$-spheres}, which for $k\geq 5$ vanishes if and only if the $k$-sphere $S^k$ has a unique differentiable structure, courtesy of the $h$-cobordism theorem of Smale \cite{smale_hcobordism}. These groups are also related to $\pi_\ast\Sph$ via the isomorphism $\Theta_{4k} \simeq \pi_{4k} \Sph/\sim$ and the short exact sequence
\[0 \to \Theta_{2k+1}^\bp \to \Theta_{2k+1} \to \pi_{2k+1} \Sph/\sim \to 0,\]
where $\Theta_{2k+1}^\bp$ is the group of homotopy spheres which bound a paralellisable manifold, and $\pi_\ast\Sph/\sim$ is the quotient of the stable stems by the image of the $J$-homomorphism, a subgroup of the $v_1$-periodic families of Adams. As $v_2$-periodic families survive this quotient, the next result follows from \Cref{mainthm:nonzerofamilies}, building upon similar $v_2$-periodic families found in \cite[Cor.1.3]{exoticwithcokerJ}, \cite[Cor.1.5]{hurewicztmf}, and \cite{bbq,bq}.

\begin{coralph}\label{maincor:exotic}
    There are exotic spheres in all dimensions congruent to $72$, $144$, and $168$ modulo $192$. 
    There are very exotic spheres, meaning exotic spheres which do not bound a parallelisable manifold, in all dimensions congruent to $143$, $145$, and $169$ modulo $192$.
\end{coralph}

In particular, the existence of exotic spheres in these three new congruence classes of dimensions modulo $192$ edges us closer to the conjecture that the only dimensions $d$ where $S^d$ has a unique smooth structure are $d=5,6,12,56$, and $61$; see \cite[Conj.2.2]{iwxpisph}.

\subsection{Related work}\label{ssec:relatedwork}
As mentioned above, some of the families in \Cref{tab:classes} recover previously studied families, as indicated by the Notes-column. In more detail:

\begin{remark}
    In a preliminary report 
    \cite[\textsection15]{tmfbook}, Hopkins--Mahowald mention a number of potential $v_2^{32}$-periodic families in $\pi_\ast\Sph_2$, some of which are outside the Hurewicz image of $\TMF$. This includes families in degrees congruent to $23$, $26$, $47$, and $48$ modulo 2, see Propositions 11.3, 11.9, 11.6, and 11.6 of \emph{ibid}, respectively. In the latter two cases, it is not clear what $v_2^?$-periodicity is considered. They also identify some of the generators in the stable stems in degrees $71-74$ in Propositions 11.14-15, although nothing about the periodicity of these classes is claimed. Our \Cref{mainthm:nonzerofamilies} proves and refines these claims, and also shows that these classes also generate $v_2^{32}$-periodic families.
\end{remark}

\begin{remark}
    In \cite{exoticwithcokerJ}, Behrens--Hill--Hopkins--Mahowald discuss nonzero elements in $\pi_\ast \Sph$, many of which are also shown to be $v_2^{32}$-periodic. In particular, the classes in rows $48$ and $120$ already appear in \cite[Tab.1]{exoticwithcokerJ}.\footnote{Curiously, the proof in \cite[Pr.12.1]{exoticwithcokerJ} that $\kappabar^6\neq0$ also uses one of the main characters of this article, $\J_0(3)$, there called $F(3)$.} The classes in row 144 also appear on page 49 of \emph{ibid}, although only as tentative $v_2^{32}$-periodic classes, and without proof. One could then consider \Cref{mainthm:nonzerofamilies} as confirming this suggestion.
\end{remark}

\begin{remark}\label{rmk:bbq}
    In \cite{bbq}, Bhattacharya--Bobkova--Quigley use the type 2 complex $A_1$ to produce $2$-torsion $v_2^{32}$-periodic families in the congruence classes $23$, $47$, $71$, $74$, $95$, $119$, and $167$. In \cite{bq}, Bobkova--Quigley exhaust this technique to produce $\eta$-torsion $v_2^{32}$-periodic families in the congruence classes $73$ and $120$. Whereas our families are produced explicitly from a generator in $\pi_*\Sph$ and a self-map, the methods of these authors show the existence of a family without identifying it explicitly. Nonetheless, Bhattacharya--Bobkova--Quigley speculate in \cite[Rmk.1.1]{bbq} on the generators of each of their families. If these speculations hold, then each of their families contains a common subfamily with one of the families of \Cref{mainthm:nonzerofamilies} by \Cref{cortoinfiniteagreement}. Similarly, we speculate that each of the Bobkova--Quigley families contains a common subfamily with one of our families in congruence classes $73$ and $120$.
\end{remark}

In proving \Cref{mainthm:nonzerofamilies}, we also refine the above works by first identifying lifts of these nonzero $v_2^{32}$-periodic families in the \emph{synthetic Hurewicz image} of $\TMF$, before proving that they also survive the ANSS for $\Sph$.

\subsection{Proof of \Cref{mainthm:nonzerofamilies}}\label{ssec:proofsection}
Our proof of \Cref{mainthm:nonzerofamilies} proceeds in two main steps. First, we compute some of the image of the ANSS of $\Sph$ in the ANSS for $\TMF$, the latter of which is well-known, although still complicated; see \cite{smfcomputation} for more background. We formulate this with Pstragowski's category of \emph{$\MU$-synthetic spectra} \cite{syntheticspectra}, a convenient, and often necessary, homotopy-theoretic category of spectral sequences. To determine families in this \emph{synthetic Hurewicz image} of $\TMF$, we combine information from the classical Hurewicz image \cite{hurewicztmf} of $\TMF$ with the use of \emph{modified ASSs}. In particular, we identify a number of families in the synthetic Hurewicz image of $\TMF$ not previously known; see \Cref{lm:47,lm:71,lm:95}. This gives us a huge collection of candidate classes in the ANSS of $\Sph$ whose behaviour in that of $\TMF$ we understand. As we are interested in families which vanish in $\pi_\ast\TMF$, this means that these candidate classes are all hit by differentials in the ANSS for $\TMF$.

Our second step, is to rule out the possibility that these differentials in the ANSS for $\TMF$ lift to differentials in the ANSS for $\Sph$. This would imply that these classes survive the ANSS for $\Sph$ and witness $v_2^{32}$-periodic families in $\pi_\ast\Sph$ at the prime $2$. Similar to Adams' use of Adams operations on $\KO$ to detect $v_1$-periodic families, we use operations on $\TMF$ to rule out the lifts of these differentials. More precisely, just as modular forms of level $N$ come with \emph{Atkin--Lehner involutions}, topological modular forms of level $N$ -- denoted as $\TMF_0(N)$ -- comes with a homotopical refinement of this involution $w\colon \TMF_0(N) \to \TMF_0(N)$, courtesy of the second-named author \cite{heckeontmf}. We then show the canonical map $p\colon \TMF \to \TMF_0(3)$ and its twist by $w$ equalises the target of these differentials, but not the source. This implies that the target lifts to the sphere, but the differential killing it does not. Most of the classes of \Cref{tab:classes} are then shown to be nonzero as they have nonzero image in this equaliser, denoted by $\J_0(3)$, reflected by a $\checkmark$ in the $\J_0(3)$-column; see \Cref{hurewiczforjzerothree}. Those remaining few classes are obtained with a filtration argument; see \Cref{thm:todabracketdetection}.

These same techniques can also be used to recover Adams' $v_1$-periodic families, as demonstrated in \cite{syntheticj}. Our techniques also suggest further avenues of research; see \Cref{ssec:questions}.

\subsection*{Outline}
In \Cref{ssec:periodicfamilies}, we recall the notation of $v_h^k$-periodic families and show that different families generated by the same class contain a common subfamily (\Cref{thm:infiniteagreement}). In \Cref{sec:synthetic}, we introduce our deleting differentials technique in synthetic spectra (\Cref{generaldifferentialkilling}) and introduce the synthetic versions of the detection spectra $\J_0(N)$ (\Cref{ssec:heighttwosyntheticstuff}). In \Cref{sec:synthetichurewicz}, we compute as much of the synthetic Hurewicz image of $\TMF_\BP$ that we need (\Cref{synthetichurewiczintext}). This is aided by Adams periodicity in $\F_2$-synthetic spectra (\Cref{sssec:adamsperiodicity}) and case-by-case arguments to lift classes to top-cells of Moore spectra (\Cref{sssec:liftsversionone,sssec:liftsversiontwo}). In \Cref{sec:periodicfamilies}, we apply the deleting differentials technique (\Cref{hurewiczforjzerothree}) and a filtration argument (\Cref{thm:todabracketdetection}) to finish the proof of \Cref{mainthm:nonzerofamilies}. In \Cref{ssec:questions}, we pose a series of questions stemming from the techniques used here.

\subsection*{Notation and terminology}

The basics of $\TMF$ are assumed, in particular its ANSS, or equivalently by \cite[Th.C]{osyn}, its descent spectral sequence (DSS). At the prime $2$, this is the localisation of the DSS for $\Tmf$ and the ANSS for $\tmf$ at $\Delta^8$, both of which are computed in \cite{smfcomputation}; helpful pictures can be found in \cite[\textsection8]{bauer}, \cite[Fig.26]{konter}, \cite[\textsection12]{tmfbook}, and \cite[Figs.A3-6]{smfcomputation}. Both Adams and Adams--Novikov spectral sequences for $\Sph$ will also be used, in which case we refer to the charts of Isaksen--Wang--Xu \cite{ANcharts,ASScharts}. The reader may want all three spectral sequences in front of them when reading \Cref{sec:synthetichurewicz,sec:periodicfamilies}.

\subsection*{Acknowledgements}
Thank you to Prasit Bhattacharya for enlightening conversations at the beginning of this project and discussing his work with Bobkova and Quigley \cite{bbq}. Thank you to Prasit and Lennart Meier for their comments on a draft too. Also thank you to Gabriel Angelini-Knoll and William Balderrama for helpful discussions and suggestions, and Sven van Nigtevecht for ongoing collaborations regarding synthetic modular forms.

The first author was supported by NSF grant \texttt{DMS-2401918} as well as the NWO grant \texttt{VI.Vidi.193.111}. The second author is an associate member of the Hausdorff Center for Mathematics at the University of Bonn (\texttt{DFG GZ 2047/1}, project ID \texttt{390685813}). We would also like to thank the Isaac Newton Institute for Mathematical Sciences, Cambridge, for support and hospitality during the programme \emph{Equivariant homotopy theory in context} where work on this paper was undertaken. This work was supported by EPSRC grant no EP/K032208/1.

\section{Generalities on periodic families}\label{ssec:periodicfamilies}
First, let us define periodic families and discuss the relationship between periodic families generated by the same element. This follows \cite[\textsection3]{exoticwithcokerJ}.

\begin{mydef}\label{def:periodicfamilies}
    Fix a prime $p$ and a nonnegative integer $h\geq 0$.
\begin{itemize}
    \item An element $x \in \pi_s \Sph$ is \emph{$v_h$-periodic} if there exists a type $h$ finite spectrum $M$, a map $\partial:M\to\Sph^d$ for some integer $d$, and a lift $\tilde{x}$ along the map
    \[\pi_{s+d}M\to\pi_s\Sph\]
    with the property that $\tilde{x}$ has nonzero image in $v_h^{-1}\pi_s M$.
    In this case, the complex $M$ admits $v_h^k$-self map $v$ on $F$ such that $v^t(\widetilde{x}) \neq 0$ for all $t\geq 1$. We say that $x$ is \emph{$v_h^k$-periodic}, emphasising the $k$.
    \item If $x$ is $v_h^k$-periodic, then fixing $M\to \Sph^d$, the choice of $v_h^k$-self map $v$, and a choice of lift $\widetilde{x}$, we write $v^t(x) \in \pi_{s+t|v|}\Sph$ for the projection $\partial(v^t(\widetilde{x}))$. We call the sequence $(v^t({x}))_{t\geq 0}$ the \emph{$v_h^k$-periodic family generated by $x$}.
    \item More generally, given a sequence of elements $\gX=(x_0, x_1, x_2, \ldots) \subseteq \pi_\ast X$, we say that $\gX$ is a \emph{$v_h^k$-periodic family} if it is a $v_h^k$-periodic family generated by $x_0$ in the sense above. If each $x_i\neq 0$, then we say this family is \emph{nonzero}.
\end{itemize}
\end{mydef}

Hidden in the notation $(v^t(x))_{t\geq 0}$ are the choices of type $h$ complex, map $\partial$, self map $v$, and lift $\tilde{x}$. For example, we have a choice of lifting a class such as $\kappa\kappabar$ using lifts of either $\kappa$ or of $\kappabar$, as found in \cite[Lm.7.17]{hurewicztmf}, and there is no immediate reason for the associated families to agree. The periodicity theorem and consequent asymptotic uniqueness of $v_h$-self maps, see \cite[Th.9]{hopkinssmith} or \cite[Th.1.5.4]{orangebook}, shows that two $v_h^k$-periodic families with the same generator have infinite intersection.

\begin{theorem}\label{thm:infiniteagreement}
    Let $\gX=(x_1,x_2,\ldots) \subseteq \pi_\ast\Sph$ be a $v_h^s$-periodic family and $\gY=(y_1,y_2,\ldots)\subseteq \pi_\ast\Sph$ be a $v_h^t$-periodic family. If $x_1=y_1$, then there are integers $i,j$ such that for all integers $n\geq 1$ we have $x_{in} = y_{jn}$ in $\pi_\ast\Sph$.
\end{theorem}

\begin{proof}
    By definition, we have maps $\partial:M\to S^d$ and $\partial:M'\to S^{d'}$ for type $h$ complexes $M$ and $M'$, a $v_h^s$-self map $v$ of $M$ and a $v_h^t$-self map $w$ of $M'$, and classes $\tilde{x_1}\in\pi_*M$ and $\tilde{y_1}\in\pi_*M'$ with the property that $\partial(\tilde{x_1})=x_1$ and $\partial'(\tilde{y_1})=y_1$. Also by definition, $x_i=\partial(v^i\tilde{x_1})$ and $y_i=\partial'(w^i\tilde{y_1})$. Suppressing shifts, we may form the commutative diagram
    \[
    \begin{tikzcd}
        M\oplus M'\arrow[d,"\pi_1"]\arrow[r,"\pi_2"]&M\arrow[d,"\partial"]\\
        M'\arrow[r,"\partial'"]&\Sph
    \end{tikzcd}
    \]
    The type $h$ complex $M\oplus M'$ admits the $v_h$-self maps $v\oplus\id$ and $\id\oplus w$, and the uniqueness of self maps, see \cite[Lm.6.1.3]{orangebook}, guarantees that there exists positive integers $i$ and $j$ such that $(v\oplus\id)^i=(\id\oplus w)^j$ as endomorphisms of $M\oplus M'$. Therefore, for all $n\ge0$
    \begin{align*}
        x_{in}&=\partial(v^{in}x_1)\\
        &=\partial(\pi_2((v\oplus\id)^{in}(x_1,y_1))\\
        &=\partial'(\pi_1((v\oplus\id)^{in}(x_1,y_1))\\
        &=\partial'(\pi_1((\id\oplus w)^{jn}(x_1,y_1))\\
        &=\partial'(w^{jn}y_1)\\
        &=y_{jn}\qedhere
    \end{align*}
\end{proof}

\begin{cor}\label{cortoinfiniteagreement}
    For any two $v_2^{32}$-periodic families $\gX$ and $\gY$ of \Cref{tab:classes} whose generators are equal the intersection $\gX \cap \gY$ contains a $v_2^{32s}$-periodic family for some $s\geq 1$.
\end{cor}

\section{Tools from synthetic spectra}\label{sec:synthetic}

Let $E$ be a homology theory of Adams type and $\Syn_E$ denote Pstragowski's $\infty$-category of \emph{$E$-based synthetic spectra}; see \cite[Df.4.1]{syntheticspectra}. The unit of $\Syn_E$ will be written as $\Sph_E$. We use \emph{stem--filtration} grading for synthetic spectra, meaning $\pi_{s,f}$ corresponds to a $(s,f)$-location in an $E$-based Adams chart. More formally, we write $\Sigma^{s,f} \Sph_E = \Sigma^{-f}\nu\Sph^{s+f}$. In particular, the element $\tau$ lives in $\pi_{0,-1} \Sph_E$ and the $\infty$-categorical suspension has bidegree $(1,-1)$. This notation follows \cite{syntheticj,smfcomputation} and differs from \cite{syntheticspectra}.

There are two functors relating spectra to synthetic spectra
\[\Sp \xrightarrow{\nu} \Syn_E \xrightarrow{\tau^{-1}} \Sp,\]
the \emph{synthetic analogue functor} and \emph{$\tau$-inversion}, whose composite is equivalent to the identity on $\Sp$. A synthetic spectrum $X$ is a \emph{synthetic lift} of a spectrum $Y$ if there is an equivalence $Y\simeq \tau^{-1}X$. There is also a lax monoidal functor
\[\sigma\colon \Syn_E \to \Fun(\Z^\op,\Sp) = \Fil(\Sp)\]
to the $\infty$-category of \emph{filtered spectra}, which sends a synthetic spectrum $X$ to its $\tau$-tower
\[\cdots \xrightarrow{\tau} \map(\Sigma^{0,1} \Sph_E, X) \xrightarrow{\tau} \map(\Sph_E, X) \xrightarrow{\tau} \map(\Sigma^{0,-1} \Sph_E, X) \xrightarrow{\tau} \cdots;\]
lax monoidality is shown in \cite[Th.1.4]{syntheticj}, and is also discussed in \cite[\textsection1.2]{osyn} and \cite[Not.2.5]{smfcomputation}. Associated to each filtered spectrum is a spectral sequence, see \cite[\textsection1.2.2]{ha} or \cite[\textsection II.1]{hedenlund_phd}. For a synthetic spectrum $X$, the filtered spectrum $\sigma(X)$ and its associated spectral sequence will be called the \emph{signature} of $X$. In \cite[Pr.1.25]{osyn}, it is shown that the signature of $\nu X$ is the $E$-based ASS for $X$, assuming that $X$ is $E$-nilpotent complete. We will often implicitly use this fact when discussing $\sigma(\nu(X))$ and its signature.

\subsection{Deleting differentials in the signature spectral sequence}\label{ssec:removingadiff} 
To illustrate the technique of deleting differentials, we keep a basic example in mind for the reader throughout this section. While $\eta^3=0$ in $\pi_*\KO$, $\eta^3\neq0$ in $\pi_{*,*}\KO_\BP$; that is, $\eta^3$ is not in the Hurewicz image of $\KO$, but it is in the synthetic Hurewicz image of $\KO$. We will use Adams operations on $\KO$ to ``delete'' the differential killing $\eta^3$ in $\KO$ to prove that $\eta^3\neq0\in\pi_*\Sph$. Recall that the $E_2$-page of the ANSS for $\KO$ is given by $\Z[u^{\pm 2},\eta]/(2\eta)$ where $|\eta|=(1,1)$ and $|u^2|=(4,0)$. This spectral sequence is determined by the differential $d_3(u^2)=\eta^3$; these facts are shown in a purely synthetic manner in \cite[\textsection4.1]{syntheticj}.

Consider now the Adams operation $\psi^3$ on $\KO$, implicitly $2$-completed. This acts on the $E_2$-page of the ANSS for $\KO$ via the formul{\ae} $\psi^3(\eta)=\eta$ and $\psi^3(u^2)=3^2u^2$. In particular, the \emph{source} of the differential $d_3(u^2)=\eta^3$ is acted on nontrivially by $\psi^3$, but the \emph{target} is acted upon trivially. This means that the source of the differential does \textbf{not} lift to the fibre of $\psi^3-1$ as a map of spectral sequences, but the target does; see \Cref{fig:removingdiff}. The class $\eta^3$ is then detected in a (modified) ANSS for the fibre of $\psi^3-1$, often called the \emph{image-of-$J$ spectrum}, and is \textbf{not} hit by a $d_3$-differential.

\begin{figure}[h]
    \centering
    \begin{sseqpage}[name = koexample1 , y range = {0}{4}, x range = {0}{8},axes type = frame, grid = go, yscale = 0.4, xscale = 0.4
]\end{sseqpage} \qquad \qquad \qquad 
    \begin{sseqpage}[name = koexample2 , y range = {0}{4}, x range = {0}{8},axes type = frame, grid = go, yscale = 0.4, xscale = 0.4
]\end{sseqpage}
    \caption{ANSS for $\KO$ on the left and the signature of $\fib(\nu\psi^3-1\colon \nu\KO \to \nu\KO)$ on the right. The blue classes are lifts from the ANSS for $\KO$, the red classes lie in the image of the boundary map, and the group in bidegree $(3,1)$ is $\Z/8\Z$; also see \cite[Fig.12]{syntheticj}.}
    \label{fig:removingdiff}
\end{figure}

This argument then shows that $\eta^3\neq0\in\pi_* \fib(\psi^3-1)$, and thus $\eta^3\neq0\in\pi_*\Sph$. We state now an almost tautologically general case of deleted differentials. In the following proposition, one recovers the example above by taking $X=Y=\KO$, $f=\psi^3-1$, and $b=\eta^3$.

\begin{prop}\label{generaldifferentialkilling}
    Let $F\xrightarrow{i} X \xrightarrow{f} Y$ be a fibre sequence of synthetic spectra and $r\geq 2$. Fix an element $b\in \pi_{\ast,\ast}X$ with $f(b)=0$, write $\bar{b}$ for its projection modulo $\tau$, and suppose that $\bar{b}\neq 0$, so that in particular, $\bar{b}$ is a nonzero permanent cycle in $\sigma(X)$. Suppose that for all $a\in \pi_{\ast,\ast} X/\tau^{r-1}$ such that $d_r(\bar{a}) = \bar{b}$, where $\bar{a}$ is the mod $\tau$-projection of $a$, one has $f(a)\neq0$. Then for all lifts $c$ of $b$ to $F$, we have $\tau^{r-1}c\neq0$.
\end{prop}

\begin{proof}
    Suppose $\tau^{r-1}c=0$, then there must be a differential $d_{r'}(x)=\bar{c}$ for some $r'\le r$ and $x \in \pi_{\ast,\ast} F/\tau$. Since $\tau^{r'-1}b\neq0$ for $r'<r$, we must have $r'=r$, as $c$ is a lift of $b$. Then by the truncated omnibus theorem \cite[Th.2.21(1)]{smfcomputation}, it follows that there exists a lift $\widetilde{x}$ of $x$ along the map $ \pi_{\ast,\ast} F/\tau^{r-1}\to \pi_{\ast,\ast} F/\tau$. If $y$ denotes the mod $\tau$-reduction of $i(\widetilde{x})$, then $d_s(y) = i(d_s(x))$ for all $s\geq 2$. In particular, $y$ is a $d_s$-cycle for $s\leq r-1$ and $d_r(y)=i(\bar{c})=\bar{b}$. Our assumptions then imply that $f(i(\widetilde{x}))\neq 0$, which contradicts the fact that $fi=0$ as a fibre sequence.
\end{proof}

To apply this general phenomenon to detection statements, we recall the notion of the \emph{Hurewicz image} of a ring spectrum. Classically, given a ring spectrum (more generally an $\E_0$-algebra) $A$, the \emph{Hurewicz image} of $A$ is the image of $\pi_\ast \Sph$ in $\pi_\ast A$. We now introduce a synthetic variant.

\begin{mydef}\label{df:synthetichurewicz}
    For an $\E_0$-algebra $\Sph_E\to A$ in synthetic spectra, we say the \emph{synthetic Hurewicz image} of $A$ is the image of the map $\pi_{*,*}\Sph_E\to\pi_{*,*}A$.
\end{mydef}


We study now a specific application of \Cref{generaldifferentialkilling} that shows how some knowledge about synthetic Hurewicz images can be used in tandem with $\sigma$-SSs to deduce facts about classical Hurewicz images. Here we follow \Cref{generaldifferentialkilling} and write $\bar{x}$ for the mod $\tau$-reduction of an element $x$ in bigraded homotopy groups of a synthetic spectrum. The $\eta^3$ example discussed above is also an instance of the following corollary, where one takes $A=B=\KO$, $\phi=\id$, $\psi=\psi^3$, and $x=\eta^3$.

\begin{cor}\label{cor:synthetichurewiczimpractice}
    Fix synthetic $\E_0$-rings $A,B$, a pair of maps of synthetic $\E_0$-rings $\varphi, \psi \colon A \to B$, an integer $r\geq 2$, and write $F$ for the equaliser of $\varphi$ and $\psi$. Suppose we have a class $x \in \pi_{s,f} A$ such that $\tau^n x$ lies in the synthetic Hurewicz image of $A$ for some $n$, that $\varphi(x) = \psi(x)$, that $\bar{x}$ is nonzero, and that for each $a \in \pi_{s+1,f-r} A/\tau^{r-1}$ with $d_r(\bar{a}) = \bar{x}$, we have $\varphi(a) \neq \psi(a)$. Then, for any lift $y$ of $x$ to $\pi_{s,f} F$, the class $y$ is not $\tau^{r-1}$-torsion. In particular, if $\pi_{s+1,f-r-i} F/\tau = 0$ for $i\geq 1$, then the $\tau$-inversion of any element $z\in \pi_{s,f}\Sph_E$ detected by $x$ lies in the classical ($E$-nilpotent complete) Hurewicz image of the underlying $\E_0$-algebra $\tau^{-1}F$.
\end{cor}

The condition that $\varphi(x) = \psi(x)$ is often satisfied. For example, if the integer $n$ above can be taken to be zero, or if multiplication by $\tau^n$ is injective on $\pi_{s,f}A$.

\begin{proof}
    We apply \Cref{generaldifferentialkilling} to this situation by setting $F=F$, $X=A$, $Y=B$, $f=\varphi - \psi$, and $b=x$; all of the hypotheses are satisfied by assumption. This tells us that all lifts $\bar{y}$ are not $d_r$-boundaries, proving the first statement. If the groups $\pi_{s+1,f-r-i} F/\tau$ vanish, then there are no possible sources for $d_s$-differentials with possible target $\bar{y}$ in $\sigma(F)$ for $s\geq r+1$. In particular, we see that each $\bar{y}$ is not just a permanent cycle in $\sigma(F)$, but that the element $y \in \pi_{s,f} F$ is not $\tau$-power torsion. As $x \in \pi_{\ast,\ast}X$ lies in the synthetic Hurewicz image, we see that at least one of the lifts $y\in \pi_{\ast,\ast}F$ also lies in the synthetic Hurewicz image. In particular, $\tau^n y =u(z)$ for some $z\in \pi_{s,f-n} \Sph_E$, where $u\colon \Sph_E \to F$ is the unit map. As $y$ is $\tau$-torsion free, we see that $z$ is also $\tau$-torsion free, and that the $\tau$-inversion of $z$ in $\pi_\ast \Sph_E^\wedge$ is sent to the $\tau$-inversion of $y$ in $\pi_\ast \tau^{-1}F$, as desired.
\end{proof}

\begin{example}
    Aside from the above example of $\eta^3$ and $\KO$, one can also apply the above to show that at the prime $3$, the product $\al_1\be_1^2\neq 0$ in the fibre of $\psi^2-1$ acting on $\tmf_3$, hence also in the $3$-complete sphere. This is implicitly done in \cite{heighttwojat3}. In forthcoming work, the second-named author will further explore iterated fibres of Adams operations and Hecke operators acting on $\TMF$ at the prime $3$.
\end{example}

\subsection{A collection of $\BP$-synthetic spectra of height $2$}\label{ssec:heighttwosyntheticstuff}
The previous section will be applied to various $\BP$-synthetic $\E_\infty$-rings related to topological modular forms. For some background, we refer the reader to \cite{tmfbook} or \cite{handbooktmf} for a general background, and \cite{smfcomputation} for more explicit details on the homotopy groups of $\TMF$.

For this whole subsection, we will implicitly complete at a prime $p$, fix an integer $N\geq 2$ which is not divisible by $p$, and work solely in $\BP$-synthetic spectra $\Syn_\BP$.

\subsubsection{Synthetic topological modular forms}

\begin{mydef}\label{def:tmf}
    For any integer $N\geq 2$, define the synthetic $\E_\infty$-rings\footnote{Many of the arguments in this article work equally as well for $\nu\tmf$, except we do not know of a reasonable connective model for the Atkin--Lehner involution $w\colon \TMF_0(3) \to \TMF_0(3)$ of \Cref{def:jspectra}.}
    \[\TMF_\BP = \nu \TMF, \qquad \TMF_0(N)_\BP = \nu \TMF_0(N).\]
\end{mydef}

We prefer the decoration $(-)_\BP$, as also used in \cite{syntheticj}, to $\nu(-)$, as $\nu$ will also often be used for the second stable Hopf element generating $\pi_3 \Sph$. This notation also highlights that these synthetic spectra a simply \emph{a} prefered synthetic lift and plays down the role of \emph{the} classical ANSS.

\begin{remark}\label{dssvsanss}
By \cite[Ths.B-C]{osyn}, the signatures of $\TMF_\BP$ and $\TMF_0(N)_\BP$ are both the ANSSs for their respective spectra as well as the associated \emph{descent spectral sequences} (DSSs).    
\end{remark}

\begin{prop}\label{pr:liftsfortmf}
    The synthetic $\E_\infty$-rings $\TMF_\BP$ and $\TMF_0(N)_\BP$ are $\tau$-complete synthetic lifts of the $\E_\infty$-rings $\TMF$ and $\TMF_0(N)$, respectively.
\end{prop}

\begin{proof}
    This follows from \cite[Pr.A.13]{burkhahnseng} and the fact that $\TMF$ and $\TMF_0(N)$ are $\BP$-nilpotent complete.
\end{proof}

\subsubsection{Synthetic height 2 image-of-$J$ spectra}
In \cite[\textsection3]{adamsontmf}, we considered the fixed points of $\TMF$ with respect to the Adams operations constructed by either \cite[Df.2.1]{heckeontmf} or \cite[Df.6.16]{luriestheorem}. The main character in this article is the following variant.



By \cite[Df.2.5]{heckeontmf}, there is a \emph{stable Atkin--Lehner involution} (also called a \emph{stable Fricke involution} there), so a map of $\E_\infty$-rings $w_N = w\colon \TMF_0(N) \to \TMF_0(N)$ inducing the classical Atkin--Lehner involution $w_N$ on the descent spectral sequence for $\TMF_0(N)$.

\begin{mydef}\label{def:jspectra}
    Let $\J_0(N)$ be the $\E_\infty$-ring given by the limit of the diagram of $\E_\infty$-rings
    \[\J_0^\bullet(N) = \begin{tikzcd}
        {\TMF}\ar[r, shift left = 1, "{q}"]\ar[r, shift right = 1, "{p}", swap]    &   {\TMF_0(N),}
    \end{tikzcd}\]
    where $p\colon \TMF \to \TMF_0(N)$ is the canonical $\TMF$-algebra map, and $q = wp$. Equivalently, the $\E_0$-ring $\J_0(N)$ is the fibre of $q-p\colon \TMF \to \TMF_0(N)$. 
\end{mydef}

The $\E_\infty$-rings $\J_0(N)$ have previously been used in detection arguments, such as in \cite[Pr.12.1]{exoticwithcokerJ} to show that $\kappabar^6\neq0$. This fact will be reproven in \Cref{hurewiczforjzerothree}.

\begin{mydef}\label{def:syntheticjspectra}
Let $\J_0(N)_\BP$ be the synthetic $\E_\infty$-ring given by the limit of the diagram of synthetic $\E_\infty$-rings
    \[\J_0^\bullet(N)_\BP = \begin{tikzcd}
        {\TMF_\BP}\ar[r, shift left = 1, "{\nu q}"]\ar[r, shift right = 1, "{\nu p}", swap]    &   {\TMF_0(N)_\BP}
    \end{tikzcd},\]
    or equivalently as the synthetic $\E_0$-ring $\fib(\nu q-\nu p\colon \TMF_\BP \to \TMF_0(N)_\BP)$. 
\end{mydef}

\begin{prop}
    The synthetic $\E_\infty$-ring $\J_0(N)_\BP$ is a $\tau$-complete lift of $\J_0(N)$.
\end{prop}

\begin{proof}
    By \Cref{pr:liftsfortmf} and the fact that $\tau$-completeness is closed under limits shows $\J_0(N)_\BP$ is $\tau$-complete. It is clear that $\J_0(N)_\BP$ is a lift of $\J_0(N)$ as the synthetic analogue functor is a section of the $\tau$-inversion functor, and $\tau$-inversion is exact.
\end{proof}

\subsubsection{Simple computations in bigraded homotopy groups}\label{ssec:pqcomputations}
One can compute the bigraded homotopy groups of $\J_0(3)_\BP$ in any particular degree from the long exact sequence induced by its defining fibre sequence, however, these homotopy groups are very busy. As we are more interested in $\J_0(3)_\BP$ as a detection spectrum, we will only compute enough of its homotopy groups to apply \Cref{cor:synthetichurewiczimpractice}.

By identifying $\pi_{\ast,\ast}\TMF_\BP/\tau$ with the $E_2$-page of the DSS for $\TMF$, see \Cref{dssvsanss}, write $\Delta\in \pi_{24,0} \TMF_\BP/\tau$ for the usual discriminant modular form.

\begin{prop}\label{pr:injectivityofpandq}
    For $0\leq i\leq 2$ and $k\geq 1$, the element $\eta^i \Delta^k \in \pi_{i+24k,i} \TMF_\BP/\tau^4$ has nonzero image under the map $q-p\colon \colon \TMF_\BP/\tau^4 \to \TMF_0(N)_\BP/\tau^4$.
\end{prop}

To apply \Cref{cor:synthetichurewiczimpractice} in this situation, will be need to delete various $d_5$-differentials, hence why we demand on working modulo $\tau^4$ in this proposition.

\begin{proof}

    First, we refer to \cite[Pr.8.1]{levelonethree}, which states that as a map between $E_2$-pages we have
    \[(q-p)(\Delta^{2^r(2s+1)}) = a_1^{3\cdot 2^{r+1}} a_3^{2^{r+1}(4s+1)} + O(a_1^{3\cdot 2^{r+1}+1})\]
    modulo $2$. These elements are all nonzero on the $E_2$-page of $\sigma(\TMF_0(3)_\BP)$. As the reduction maps $\TMF_\BP/\tau^4 \to \TMF_\BP/\tau$ and $\TMF_0(3)_\BP/\tau^4\to \TMF_0(3)_\BP/\tau$ are both isomorphisms in degrees $\pi_{24k,0}$, this proves our statement for $i=0$. For $i\neq0$, we use the fact that $\eta^i$ comes from $\nu\Sph$ and $q-p$ is $\nu\Sph$-linear, which leads to $(q-p)(\eta^i \Delta^k) = \eta^i (q-p)(\Delta^k)$. This expression is nonzero by inspection of $\pi_{\ast,\ast}\TMF_0(3)_\BP/\tau^4$.
\end{proof}

We will often use the following simple facts about $\J_0(N)_\BP$.

\begin{lemma}\label{lm:checkboard}
    The $\BP$-synthetic spectrum $\J_0(N)_\BP$ is \emph{even}, so in particular, $\sigma(\J_0(N)_\BP)$ is concentrated in even total degrees and can only support odd differentials.
\end{lemma}

In other words, $\sigma(\J_0(N)_\BP)$ obeys a ``checkerboard'' pattern.

\begin{proof}
    This is true for $\TMF_\BP$ and $\TMF_0(N)_\BP$ and this condition is closed under fibres.
\end{proof}

\begin{lemma}\label{lm:nothingbelowzero}
    For $f\leq -1$, the groups $\pi_{\ast,f}\J_0(N)_\BP/\tau$ all vanish.
\end{lemma}

\begin{proof}
    This follows from this fact for $\TMF_\BP$ and $\TMF_0(N)_\BP$, which is clear as their signatures are classical ANSSs. 
\end{proof}

\begin{lemma}\label{lm:filtrationzero}
    For all $n\neq0$, the groups $\pi_{n,0}\J_0(N)_\BP/\tau$ vanish. In other words, the zero line of $\sigma(\J_0(N)_\BP)$ is zero away from the zero stem.
\end{lemma}

\begin{proof}
    From the long exact sequence on mod $\tau$-bigraded homotopy groups defining $\J_0(N)_\BP$, the zero line of $\sigma(\J_0(N)_\BP)$ comes from classes in the zero line of $\sigma(\TMF_\BP)$ which lie in the kernel of $q-p$. One of the key features of the Atkin--Lehner involution on $\TMF_0(N)$ is that $w^2\simeq \psi^N$; see \cite[Th.C]{heckeontmf}. In particular, it suffices to compute the kernel of the map
    \[(w+1)(q-p) \simeq \psi^Np - q + q - p \simeq \psi^Np-p \simeq p\psi^N - p \simeq p(\psi^N-1)\]
    on the zero line of $E_2$-pages; the second-to-last homotopy is the naturality of Adams operations from \cite[(2.2)]{heckeontmf}. The map $p$ is injective on $0$-lines between $\sigma(\TMF_\BP)$ and $\sigma(\TMF_0(N)_\BP)$, as a modular form of level $\Ga_0(1)$ is a modular form of level $\Ga_0(N)$. It suffices to show then that $\psi^N-1$ has trivial kernel on the zero line of $\sigma(\TMF_\BP)$ away from the zero stem. For each $x\in \pi_{s,0}\TMF_\BP/\tau$, we have $\psi^N(x)=N^{\tfrac{s}{2}}x$ by \cite[Cor.2.12]{heckeontmf}, and that all of these abelian groups $\pi_{s,0}\TMF_\BP/\tau$ are torsion-free. In particular, such an $x$ is in the kernel of $\psi^N-1$ if and only if it lies in the zero stem $s=0$. This finishes the proof.
\end{proof}

\section{On the synthetic Hurewicz image of topological modular forms}\label{sec:synthetichurewicz}
With all of our synthetic spectra at hand, we can begin with our computations. First, we show that all of the classes in \Cref{tab:classes} admit $\tau$-power torsion lifts to $\TMF_\BP$ and lie in the synthetic Hurewicz image; see \Cref{synthetichurewiczintext}. 

In this section, we implicitly complete at the prime $2$.
\subsection{Easy classes in the synthetic Hurewicz image}
Write $\Delta^8 \in \pi_{192,0}\TMF_\BP$ for the unique class whose mod $\tau$-reduction is the $8$th power of the discriminant modular form $\Delta$; this exists in $\TMF_\BP$ by \cite[Lm.8.1]{smfcomputation}.

\begin{notation}\label{not:elementsfromsphere}
    Let us write
    \[\eta \in \pi_{1,1} \TMF_\BP, \qquad 
    \nu \in \pi_{3,1} \TMF_\BP, \qquad 
    \eps \in \pi_{8,2} \TMF_\BP, \qquad 
    \kappa \in \pi_{14,2} \TMF_\BP,\]
    for the image of the similarly named elements from the sphere $\Sph_\BP$; these have nonzero image in $\TMF_\BP$ courtesy of \cite[\textsection5]{smfcomputation}. We also write $\kappabar \in \pi_{20,4} \TMF_\BP$ for the unique class such that $\tau^2\kappabar$ is the image of the element $g \in \pi_{20,2} \Sph_\BP$; as discussed in \cite[Not.5.9]{smfcomputation}, this element $g$ in $\Sph_\BP$ is simply a fixed choice of class in the Toda bracket $\angbr{\kappa, 2, \eta, \nu}$. Classes in $\pi_{\ast,\ast} \TMF_\BP$ are denoted by their $E_2$-page representatives\footnote{All of the classes in \Cref{tab:classes} will admit unique non-$\tau$-divisible lifts to $\TMF_\BP$, i.e., there are no contributions from higher filtration classes in these particular degrees.}, with square brackets signifying indecomposable classes whose mod $\tau$-reduction has the suggested product decomposition.
\end{notation}

\begin{theorem}\label{synthetichurewiczintext}
    Every generator of \Cref{tab:classes}, as well as all of their $\Delta^{8k}$-multiples for all $k\geq 1$, admit unique $\tau$-power torsion $\BP$-synthetic lifts to $\TMF_\BP$ of indicated AN-filtration. All of these classes, or some $\tau$-multiple, lie in the synthetic Hurewicz image of $\TMF_\BP$.
\end{theorem}

\begin{proof}
    Except for the classes in degrees congruent to $47$, $71$, and $95$ modulo $192$, all of the classes in \Cref{tab:classes} are products of $\tau$-torsion free elements which $\tau$-invert to elements in the classical Hurewicz image computed in \cite{hurewicztmf}; this is precisely as they were chosen. In particular, some $\tau$-power multiple, less than their filtration degree, lies in the synthetic Hurewicz image of $\TMF_\BP$. The only classes left to verify are those in the exceptional degrees mentioned above, which are covered in \Cref{lm:47,lm:71,lm:95}.
\end{proof}


\subsection{$\F_2$-synthetic Moore spectra and Adams periodicity}\label{sssec:adamsperiodicity}
Finishing the proof of \Cref{synthetichurewiczintext} will take some work, although most of the statements are proven by analysing the ASS charts of \cite{ASScharts} and observing that many the bidegree is question is purely $\tau$-power torsion together with information from Adams $\mu$-family or the image-of-$J$.

More precisely, we use $\Syn_{\F_2}$ to easily access modified ASSs, similar to those in \cite{mmmm,exoticwithcokerJ,hurewicztmf}, and study homotopy groups of $\F_2$-synthetic Moore spectra; see \cite[Ex.1.8]{syntheticj} for an appearance of these modified spectral sequences in the synthetic language.

\begin{mydef}\label{modifedadamsspectralsequences}
    For $i,j\geq 1$, define the $\F_2$-synthetic spectra $M(h_0^i)$ and $M(h_0^i, v_1^j)$ via the cofibre sequences in $\Syn_{\F_2}$
    \[\Sigma^{0,i}\Sph_{\F_2} \xrightarrow{h_0^i} \Sph_{\F_2} \to M(h_0^i) \xrightarrow{\partial_{0^i}} \Sigma^{1,i-1} \Sph_{\F_2}\]
    \[\Sigma^{2j,j}M(h_0^i) \xrightarrow{v_1^j} M(h_0^i) \to M(h_0^i,v_1^j) \xrightarrow{\partial_{1^j}} \Sigma^{2j+1,j-1}M(h_0^i),\]
    when the self-maps $v_1^j$ exist, for example, for $(i,j)=(1,4)$, $(2,4)$, and $(3,8)$; see \cite{mmmm,hurewicztmf}. Define the $\BP$-synthetic spectra $M(2^i)_\BP = \nu M(2^i)$ and $M(2^i,v_1^j)_\BP = \nu M(2^i,v_1^j)$ for appropriate $j$, and the map $v_2^{32}\colon \Sigma^{192,0} M(2^i,v_1^j)_\BP \to M(2^i,v_1^j)_\BP$ by applying the synthetic analogue to the known self-maps of Behrens--Hill--Hopkins--Mahowald \cite{mmmm} and Behrens--Mahowald--Quigley \cite[Th.1.8]{hurewicztmf}.
\end{mydef}

As multiplication by $2^i$ on $\Sph$ and $v_1^j$ on $M(2^i)$ induce short exact sequences on $\BP$-homology, then \cite[Lm.4.23]{syntheticspectra} states that the sequences
\[\Sph_{\BP} \xrightarrow{2^i} \Sph_{\BP} \to M(2^i)_\BP, \qquad \Sigma^{2j,0} M(2^i)_\BP \xrightarrow{v_1^j} M(2^i)_\BP \to M(2^i,v_1^j)_\BP,\]
equipped with their natural nullhomotopies, are cofibre sequences in $\Syn_\BP$.

To have a better understanding of these modified ASSs, we will need to use a periodicity operator similar to that of Adams \cite[Th.1.2]{adamsperiodicity}.

\begin{mydef}
    Let $1\leq i\leq 3$. For each $h_0^i$-torsion class $x$ in $\pi_{s,f} \Sph_{\F_2}$, we define the set $P_i(x) \subseteq \pi_{s+8,f+4} \Sph_{\F_2}$ by the formula $P(x) = \partial_{0^i}(v_1^4 \partial_{0^i}^{-1}(x))$ using the boundary maps in the diagram
    \[\begin{tikzcd}
        {\pi_{s+1,f+i-1} M(h_0^3)}\ar[r, "{\partial_{0^i}}"]\ar[d, "{v_1^4}"]   &   {\pi_{s,f} \Sph_{\F_2}} \\
        {\pi_{s+9,f+i+3} M(h_0^3)}\ar[r, "{\partial_{0^i}}"]                    &   {\pi_{s+8,f} \Sph_{\F_2},}
    \end{tikzcd}\]
    where $v_1^4\colon \Sigma^{8,4} M(h_0^i) \to M(h_0^i)$ is the self-map from \Cref{modifedadamsspectralsequences}. If $P(x)$ is a singleton set, it is automatically $h_0^i$-torsion, and we write $P_i^2(x) = P_i(P_i(x))$. For $n\geq 3$, if $P_i^{n-1}(x)$ is a singleton set, then we write $P_i^n(x) = P_i(P_i^{n-1}(x))$. When $P_i(x)$ is a singleton, we say that it has \emph{zero indeterminacy}.
\end{mydef}

Our goal in this section is to see how this periodicity operator acts on elements in $\pi_{\ast,\ast}\Sph_{\F_2}$ coming from the image-of-$J$ or Adams $\mu$-family. These are standard facts about height 1 periodic families that we express in synthetic language for our purposes; see \cite[\textsection5.3]{greenbook} for a discussion of these classes in the ANSS for $\Sph$.

\begin{mydef}
For $k\geq 0$, let $\mu_{8k+1} \in \pi_{8k+1, 4k+1} \Sph_{\F_2}\simeq \F_2$ be the unique generator in this degree, and $\mu_{8k+2} = h_1 \mu_{8k+1}$. For $k\geq 1$, we define
    \begin{enumerate}
        \item $j_{8k-5} \in \pi_{8k-5,4k-3} \Sph_{\F_2} \simeq \Z/8\Z$ as a generator.
        \item $j_{8k-1} \in \pi_{8k-1, 4k-3-\nu_2(k)} \Sph_{\F_2} \simeq \Z/2^{\nu_2(k)+4} \oplus A$ as a $\tau$-torsion free generator of order $\nu_2(k)+4$, where the elements of $A$ are $\tau$-power torsion. We also demand that $j_{15}$ and $j_{23}$ both $\tau$-invert to the classical image-of-$J$. These classes $j_{8k-1}$ are only well-defined up to $\tau$-power torsion.
        \item $j'_{8k-1} = h_0^{\nu_2(k)+1} j_{8k-1} \in \pi_{8k-1, 4k-2} \simeq \Z/8\Z$. This class is uniquely defined, independent of our choice of $j_{8k-1}$.
        \item $j_{8k} \in \pi_{8k, 4k-1} \Sph_{\F_2} \simeq \F_2$ as the generator.
        \item $j_{8k+1} = h_1 j_{8k}$.
    \end{enumerate}
     Call the set $\{\mu_{8k+1}, h_1\mu_{8k+1}\}_{k\geq 0} \subseteq \pi_{\ast,\ast} \Sph_{\F_2}$ the \emph{Adams $\mu$-family} and the set
     \[\{h_0^i j_{8k-5}, h_0^i j_{8k-1}, j_{8k}, j_{8k+1}\}_{i,k-1\geq 0} \subseteq \pi_{\ast,\ast} \Sph_{\F_2}\]
     the \emph{image-of-$J$}.
\end{mydef}

\begin{prop}\label{pr:periodicity1}
    For $1\leq i\leq 3$, $k\geq 0$, and $n\geq 1$ we have the equalities in $\pi_{\ast,\ast}\Sph_{\F_2}$
    \[P_i^n(\mu_{8k+1}) = \mu_{8(k+n)+1}, \qquad P_i^n(\mu_{8k+2}) = \mu_{8(k+n)+2}.\]
    For $1\leq i\leq 3$, $k\geq 1$, and $n\geq 1$ we have the equalities in $\pi_{\ast,\ast}\Sph_{\F_2}$
    \[P_3^n(j_{8k-5}) = j_{8(k+n)-5}, \qquad P_3^n(j'_{8k-1}) = j'_{8(k+n)-1}\]
    \[P_i^n(j_{8k}) = j_{8(k+n)}, \qquad P_i^n(j_{8k+1}) = j_{8(k+n)+1}.\]
\end{prop}

\begin{proof}
    It suffices to do these computations modulo $\tau$ for the $\mu$-family and the image-of-$J$ away from degrees congruent to $1$ modulo $8$; the indeterminacy in each case vanishes. In these cases, an explicit computation as done in \cite[Th.1.2]{adamsperiodicity} shows the desired result. The case for classes $j_{8k+1} = h_1 j_{8k}$ then follows from $j_{8k}$ by $h_1$-multiplication.
\end{proof}

For the following two propositions, we write $\overline{x} \in \pi_{s,f}M(h_0^3)$ for the projection of an element $x\in \pi_{s,f}\Sph_{\F_2}$.

\begin{prop}\label{pr:periodicity2}
    For $k\geq 0$ and $n\geq 1$ we have the equalities in $\pi_{\ast,\ast} M(h_0^3)$
    \[v_1^{4n}\overline{\mu}_{8k+1} = \overline{\mu}_{8(k+n)+1}, \qquad v_1^{4n}\overline{\mu}_{8k+2} = \overline{\mu}_{8(k+n)+2}.\]
    For $k\geq 1$ and $n\geq 1$ we have the equalities in $\pi_{\ast,\ast} M(h_0^3)$
    \[v_1^{4n}\overline{j}_{8k} = \overline{j}_{8(k+n)}, \qquad v_1^{4n}\overline{j}_{8k+1} = \overline{j}_{8(k+n)+1}, \qquad v_1^{4n}\overline{j}_{8k-5} = \overline{j}_{8(k+n)-5}\]
\end{prop}

\begin{proof}
    Let $x\in \pi_{s,f} \Sph_{\F_2}$ and write $\overline{x} \in \pi_{s,f} M(h_0^3)$ for its projection. Then we have $\partial_{0^3}(v_1^4 \overline{x}) = h_0h_3 x$ in $\pi_{s+7,f+2} \Sph_{\F_2}$. Indeed, this comes from a diagram chase and the fact that $h_0 h_3 = \partial_{0^3}(v_1^3)$; see \cite[Lm.2.6]{heighttwojat3} for the analogous statement at the prime $3$. In particular, we see that $\partial_{0^3}(v_1^4 \overline{\mu}_{8k+1}) = h_0h_3 \mu_{8k+1} = 0$, as $\mu_{8k+1}$ are all $h_0$-torsion. Consider the exact sequence
    \[\coker(\pi_{8k+9, 4k+2} \Sph_{\F_2} \xrightarrow{h_0^3} \pi_{8k+9, 4k+5} \Sph_{\F_2}) \to \pi_{8k+9, 4k+5} M(h_0^3) \xrightarrow{\partial_{0^3}}  \pi_{8k+8, 4k+3} \Sph_{\F_2}.\]
    We just saw that $v_1^4 \overline{\mu}_{8k+1}$ lies in the kernel of $\partial_{0^3}$, hence it is hit by the unique class $\mu_{8k+9}$ generating the cokernel on the left-hand side. By induction, this gives the first equality, and multiplication by $h_1$ gives us the second equality. The same arguments hold in the $j_{8k}$- and $j_{8k-5}$-cases, and the $j_{8k+1}$-case then follows by $h_1$-multiplication.
\end{proof}

\begin{prop}\label{pr:periodicity3hard}
    For $k\geq 3$ and $n\geq 1$ we have the equality in $\pi_{\ast,\ast} M(h_0^3)$
    \[\tau^a v_1^{4}\overline{j}_{8k-1} \equiv \tau^b \overline{j}_{8k+7}\]
    up to $\tau$-power torsion, where either $a=0$ or $b=0$ such that the degrees match.\footnote{This means that the first equality, up to $\tau$-power torsion, takes place in $\pi_{8k+7,f} M(h_0^3)$, where $f$ is the minimum of the filtration of $\mu_{8k+7}$ and the filtration of $\mu_{8k-1}$ plus $4$.}
\end{prop}


\begin{proof}
    Let us write $x_k \in \pi_{8k-1}\Sph$ for the $\tau$-inversion of $j_{8k-1}$. The unit map $\Sph \to \j^1$, where $\j^1$ is the fibre of $\psi^3-1\colon \ko\to \ko$, sends $x_k$ to the generator of $\pi_{8k-1} \j^1 \simeq \Z/2^{\nu_2(k)+4}\Z$ which we denote by $y_k$; see \cite{imageofstableJ} or \cite[Th.A]{syntheticj} for a synthetic proof. Moreover, the class $y_k$ in $\j^1$ hit by the class $\be^k \in \pi_{8k}\ko$ under the boundary map $\partial\colon \ko[-1] \to \j$. Tensoring this map with the $v_1^4$-self map on $M(8)$ then gives the commutative diagram
    \[\begin{tikzcd}
        {\pi_{8k} \ko/8}\ar[r, "\partial"]\ar[d, "{v_1^4}"] &   {\pi_{8k-1} \j^1/8}\ar[d, "{v_1^4}"]   \\
        {\pi_{8k+8} \ko/8}\ar[r, "\partial"]                &   {\pi_{8k-1} \j^1/8.}
    \end{tikzcd}\]
    By construction of $v_1^4$, the left-hand vertical map is multiplication by the mod $8$ reduction of the Bott class $\be$. In particular, the commutativity of this diagram states that $v_1^4\overline{y}_{k} = \overline{y}_{k+1}$. In particular, we have $v_1^4\overline{x}_k \equiv \overline{x}_{k+1}$ modulo elements in the kernel of the unit map $\Sph \to \j^1$. In $\F_2$-synthetic spectra, we then see that we must have
    \begin{equation}\label{eq:littletaupowertorsionjthingy}\tau^a v_1^{4}\overline{j}_{8k-1} \equiv \tau^b \overline{j}_{8k+7}\end{equation}
    modulo $\tau$-power torsion and elements in the kernel of $M(h_0^3) \to M(h_0^3)\otimes \j^1_{\F_2}$, where $\j^1_{\F_2}$ is the fibre of $\nu(\psi^3-1)\colon \nu\ko\to \nu\ko$ as in \cite[Df.5.3]{syntheticj}, up to a connective cover. In the exact sequence
    \[\pi_{8k+7,f} \Sph_{\F_2} \to \pi_{8k+7,f} M(h_0^3) \xrightarrow{\partial_{0^3}} \pi_{8k+6,f-2} \Sph_{\F_2},\]
    we see that both the left and right groups are $\tau$-power torsion courtesy of the assumption that $k\geq 3$, so for degree reasons we have (\ref{eq:littletaupowertorsionjthingy}) up to $\tau$-power torsion, as desired.
\end{proof}

\subsection{Lifts to the top-cell of $M(8,v_1^8)$}\label{sssec:liftsversionone}
We now come to our first indecomposable $\tau$-power torsion element in the synthetic Hurewicz image of $\TMF_\BP$.

\begin{prop}\label{lm:47}
    For all $k\geq 0$, the element $\kappabar[2\nu\Delta]\Delta^{8k}$ in $\pi_{47+192k,5}\TMF_\BP$, or a nonzero $\tau$-power multiple, lies in the synthetic Hurewicz image.
\end{prop}

To prove \Cref{lm:47}, we first show that $\kappabar[2\nu\Delta]$ lies in the synthetic Hurewicz image (\Cref{lm:47one}), then we show that the associated class in $\pi_{47} \Sph$ lifts to the top-cell of $M(8,v_1^8)$ using $\F_2$-synthetic spectra (\Cref{lm:47two}), before we finally find a $\BP$ synthetic lift of this lift to the top-cell. The proof of \Cref{lm:71} also follows this outline.

\begin{lemma}\label{lm:47one}
    The class $\kappabar[2\nu\Delta] \in \pi_{47,5}\TMF_\BP$ lies in the synthetic Hurewicz image.
\end{lemma}

\begin{proof}
    Inside $\pi_{47,5} \TMF_\BP$ we have the Toda bracket $\angbr{\tau^2\kappabar^2,\tau^2\nu, 2\nu}$ which has zero indeterminacy. Using a synthetic version of Moss' theorem, see \cite[Th.3.16]{smfcomputation}, and the $E_5$-page of the DSS for $\TMF$, we compute that our element $\kappabar[2\nu\Delta]$ is given by the above Toda bracket. Inside $\pi_{47,5} \Sph_\BP$ we have the similar Toda bracket given as $\angbr{g^2, \tau^2\nu, 2\nu}$, where $g^2\in \pi_{40,6}\Sph_\BP$ detects $\kappabar^2$. The synthetic Moss' theorem identifies a class inside this bracket as $\Delta h_2^2 e_0$ using the notation of the $E_6$-page of the ANSS for $\Sph$ of \cite{ANcharts}. In particular, this Toda bracket in $\Sph_\BP$ is nonempty, and by linearity of Toda bracket is sent to the class $\kappabar[2\nu\Delta]$ in $\TMF_\BP$ as desired as this is true for the constituents of the above bracket; see \Cref{not:elementsfromsphere}.
\end{proof}

Looking now at the $E_\infty$-page of the ANSS for $\Sph$ of \cite{ANcharts}, we see that the class $z:=\Delta h_2^2 e_0$ has the property that both $2z$ and $\eta z$ are $\tau$-torsion free. Set $x_{47} := \tau^{-1}z\in \pi_{47}\Sph$.

\begin{lemma}\label{lm:47two}
    The class $x_{47}\in \pi_{47}\Sph$ admits an $\F_2$-synthetic lift to $y_{47}\in \pi_{47,10}\Sph_{\F_2}$. Moreover, up to the image-of-$J$, $x_{47}$ admits a lift to the top-cell of $M(8,v_1^8)$.
\end{lemma}

\begin{proof}
    Let us abbreviate $x_{47}=x$ for this proof. As $x$ is nonzero, we know there is a nonzero $\tau$-torsion free class $y$ in $\pi_{47,f}\Sph_{\F_2}$ which lifts $x$ and is not divisible by $\tau$. As $x$ supports multiplication by $2$ and $\eta$, then $y$ must also support multiplication by $\tau^r h_0$ and $\tau^r h_1$ for all $r\geq 0$. By inspection of the $E_\infty$-page of \cite{ASScharts} in stem $47$, we see that $f\geq 10$, so we may $f=10$ by taking $\tau$-powers if necessary.
    
    By further inspection, we see that $y \in \pi_{47,10}\Sph_{\F_2}$ is $h_0^3$-torsion; it supports $h_0$-multiplication to $\tau^2 P\Delta h_1 d_0$ in $\pi_{47,11}\Sph_{\F_2}$, and for degree reasons, and the fact that the generator of the image-of-$J$ in filtration $20$ is not $h_0$-divisible, this class is $h_0$-torsion. Hence $y$ admits a lift to a class $y[1,2] \in \pi_{48,12}M(h_0^3)$ along the projection to the top cell $M(h_0^3)\to \Sph^{1,2}$ (the square brackets in $y[1,2]$ indicate the shift in bidegree).
    
    Consider $v_1^4 y[1,2]$ in $\pi_{56,16} M(h_0^3)$, which fits into the short exact sequence
    \begin{equation}\label{eq:firstimJcorrection} \coker(h_0^3\colon \pi_{56,13}\Sph_{\F_2} \to \pi_{56,16} \Sph_{\F_2}) \xrightarrow{\iota} \pi_{56,16} M(h_0^3) \xrightarrow{\partial_{0^3}} \left(\pi_{55,14} \Sph_{\F_2}\right)[h_0^3],\end{equation}
    where $A[h_0^3]$ indicates the $h_0^3$-torsion in $A$. Looking at the $E_2$-page of \cite{ASScharts} in bidegree $(55,14)$ (indicated with an $il$), we see that up to $\tau$-power torsion, the kernel of $h_0^3$ on $\pi_{55,14}\Sph_{\F_2}$ is spanned by $\tau^{12} P^6 (h_0h_3)$. Suppose that $\partial_{0^3}(v_1^4y[1,2]) = \tau^{12} P^6(h_0h_3)$ modulo $\tau$-power torsion. We then define $y' = y - \tau^{12} P^5(h_0h_3)$; this is an $\F_2$-synthetic lift of $x+j$ where $j$ lies in the image-of-$J$. The class $P^5(h_0h_3)$ is $h_0^3$-torsion, hence admits a lift $P^5(h_0h_3)[1,2]$ along $\partial_{0^3}$. Therefore, $y'$ admits the lift
    \[y'[1,2] = y[1,2] - \tau^{12} P^5(h_0h_3)[1,2] \in \pi_{48,12} M(h_0^3).\]
   such that
    \[\partial_{0^3}(v_1^4 y'[1,2]) = \tau^{12} P^6(h_0h_3) - \tau^{12} P^6(h_0h_3) = 0\]
    up to $\tau$-power torsion, as $\partial_{0^3}(v_1^4 P^5(h_0h_3)[1,2]) = P^6(h_0h_3)$ by \Cref{pr:periodicity1}. 
    
    In particular, we may assume that $\partial_{0^3}(v_1^4 y[1,2])$ vanishes up to $\tau$-power torsion, so that it lies in the image of $\iota$. Looking now at the $E_2$-page of \cite{ASScharts} in bidegree $(56,16)$, we see that up to $\tau$-power torsion, the image of $\iota$ in (\ref{eq:firstimJcorrection}) is spanned by $\tau^{11} \overline{P^6 c_0}$, where $P^6 c_0$ generates the image-of-$J$ in this degree and $\overline{(-)}$ indicates its projection to $M(h_0^3)$. If $v_1^4 y[1,2]$ is not $\tau$-power torsion, then up to $\tau$-power torsion it must be $\tau^{11} \overline{P^6 c_0}$. This time, we can alter $y[1,2]$ by setting $y[1,2]' = y[1,2] - \tau^{11} \overline{P^5 c_0}$, and $\partial_{0^3}(y[1,2]')=\partial_{0^3}(y[1,2])=y$. Courtesy of \Cref{pr:periodicity2}, we see that $v_1^4 \overline{P^5 c_0} = \overline{P^6 c_0}$, hence
    \[v_1^4 y[1,2]' = \tau^{11} \overline{P^6 c_0} - \tau^{11} \overline{P^6 c_0} = 0\]
    up to $\tau$-power torsion. In particular, the $\tau$-inversion of $y[1,2]'$, written as $x[1] \in \pi_{48} M(8)$ is $v_1^4$-torsion, hence it admits a lift to $x[18]$ in $\pi_{65} M(8,v_1^8)$.
\end{proof}

\begin{proof}[Proof of \Cref{lm:47}]
    Finally, we lift the class $z\in\pi_{47,5}\Sph_\BP$ of \Cref{lm:47one} to $M(8,v_1^8)_\BP$ by producing an appropriate lift of $x[18]$ to a $\BP$-synthetic class. The element $x[18]$ is detected by a non-$\tau$-divisible class $\widetilde{w}\in \pi_{65,f} M(8,v_1^8)_\BP$ such that the element $w:=\partial_{0^3} \partial_{1^8} (\widetilde{w})\in \pi_{47,f+2} \Sph_\BP$ is a lift of $x$ along $\tau$-inversion, so that $\tau^aw=\tau^bz$ modulo $\tau$-torsion for some $a,b\ge0$. The charts of \cite{ANcharts} show that the bidegrees $\pi_{47,>5}\Sph_\BP$ consist entirely of $\tau$-torsion, and $\pi_{47,\le 5}\Sph_\BP$ is $\tau$-torsion free. It follows that $w=\tau^r z$ for $r\ge0$, so that $w$ is mapped to $\tau^r \kappabar[2\nu\Delta]$ in $\TMF_\BP$ for some $r\geq 0$. It follows that $0\leq f\leq 3$ since $\widetilde{w}$ is not $\tau$-divisible, so that $r<4$ and $\tau^r \kappabar[2\nu\Delta]\neq0$.
    
    We now generate a family of elements $w_k \in \pi_{47+192k,f+2} \Sph_\BP$ inductively by setting $\widetilde{w}_0:=\widetilde{w}$,
    \[\widetilde{w}_{k+1} := v_2^{32} \widetilde{w}_k \in \pi_{65+192(k+1),f} M(8,v_1^8),\]
    and finally $w_k := \partial_{0^3}\partial_{1^8}(\widetilde{w}_k)$. We claim that the elements $\tau^r\kappabar[2\nu\Delta]\Delta^{8k}$ in $\TMF_\BP$ are in the synthetic Hurewicz image and detected by $w_k$, up to a sign. We have already seen this above for $k=0$. In general, this follows from the fact that multiplication by $v_2^{32}$ on the homotopy groups of $\TMF_\BP\otimes M(8,v_1^8)$ induced by the self-map of $M(8,v_1^8)$ equivalent (up to a sign) to the $\TMF_\BP$-module map induced by multiplication by $\Delta^8$. This in turn is a consequence of the equation
    \begin{equation}\label{eq:deltacongruence}\Delta = v_2^3v_1^3 - 27v_2^4\end{equation}
    in the ring of ($2$-local) modular forms, a reflection of the fact that $\Delta = a_3^3a_1^3-27a_3^4$ in the ring of modular forms of level $\Ga_0(3)$ and the equalities $a_1=v_1$ and $a_3 = v_2$; see \cite[Lm.1]{laureskonetmf}. This finishes the proof.
\end{proof}

\begin{warn}
    It is tempting to work solely in either $\Syn_\BP$ or $\Syn_{\F_2}$ to prove \Cref{lm:47}, but this has drawbacks. For example, the class $g \in \pi_{20,4} \Sph_{\F_2}$ is $h_0^3$-torsion, hence it lifts to a class $g[1,2]\in \pi_{21,6} M(h_0^3)$. This second class is \textbf{not} $v_1^8$-torsion though---it is only true that $v_1^8 g[1,2]$ is $\tau$-power torsion, hence it lifts classically but \textbf{not} $\F_2$-synthetically. A similar problem shows that $y[1,2]\in \pi_{48,12} M(h_0^3)$ also does not lift to $M(h_0^3,v_1^8)$. The most elegant solution is presumably a \emph{bisynthetic} approach.
\end{warn}

There is an arguably more elegant alternative proof of \Cref{lm:47}; we prefer the proof above as these techniques more easily generalise to \Cref{lm:71}.

\begin{remark}
    Classically, there is a Toda bracket $\eta^2 = \angbr{2,\eta,2} \in \pi_2 \Sph$ and this also lifts to a synthetic Toda bracket $\tau^2\eta^2 = \angbr{2,\eta,2} \in \pi_{2,2} \Sph_\BP$. In particular, if a class $x\in \pi_{s,f} X$ is $2$-torsion for synthetic spectrum $X$, then $\tau^2\eta^2 x$ is divisible by $2$. Indeed, this is just the juggling formula for synthetic Toda brackets
    \[\tau^2\eta^2 x = \angbr{2,\eta,2}x = 2\angbr{\eta, 2, x};\]
    see \cite[Pr.3.13(1)]{smfcomputation}, for example. Applying this to $\eta^2\kappabar[\eta\Delta]$, whose $\tau$-power multiple lies in the synthetic Hurewicz image of $\TMF_\BP$, is $2$-divisible in both $\TMF_\BP$ and (the appropriate $\tau$-power multiple) $\Sph_\BP$. In particular, this shows that a $\tau$-power multiple of $\kappabar[2\nu\Delta]$ lies in the synthetic Hurewicz image of $\TMF_\BP$.
\end{remark}

\begin{prop}\label{lm:71}\
    For all $k\geq 0$, the element $\kappabar[\nu\Delta^2]\Delta^{8k}$ in $\pi_{71+192k,5}\TMF_\BP$, or a nonzero $\tau$-power multiple, lies in the synthetic Hurewicz image.
\end{prop}

\begin{lemma}\label{lm:71one}
    The class $\kappabar[\nu\Delta^2] \in \pi_{71,5}\TMF_\BP$ lies in the synthetic Hurewicz image.
\end{lemma}

\begin{proof}
    First, note that $\pi_{65} \TMF$ contains a class $[\nu\kappa\Delta^2]$ which lies in the Hurewicz image (\cite[Lm.7.21(2)]{hurewicztmf}) and has AN-filtration $3$. We claim that the corresponding class in $\pi_{65,3}\TMF_\BP$ admits a lift to $x\in\pi_{65,3}\Sph_\BP$ such that $\nu^2x$ is divisible by $4$, modulo elements in the kernel of the Hurewicz map to $\TMF_\BP$. This suffices to prove the proposition, as $\nu^2[\nu\kappa\Delta^2]$ also lies in the synthetic Hurewicz image of $\TMF_\BP$, and $\nu^2[\nu\kappa\Delta^2] = 4\kappabar[\nu\Delta^2]$.
    
    Let $\pi(y)$ denote the mod $\tau$ projection of a class $y$. Referring to the $E_\infty$-page of the ANSS for $\Sph$ of \cite{ANcharts}, we see that for any lift $x$ of $[\nu\kappa\Delta^2]$, we have $\pi(x)=\sum\limits_{i=1}^4\eps_ix_i$, where $x_i$ are the generators of bidegree $(65,3)$ of the $E_\infty$-page labelled from left to right, and $\eps_i=0,1$. We claim that $\eps_1=1$. Indeed, the class labelled $h_5^2$ in bidegree $(62,2)$ goes to zero in $\pi_{62,2}\TMF_\BP/\tau$ since it is a permanent cycle in filtration $2$ and $\pi_{62}\TMF=0$. It follows then that $\eps_2=0$. Now since $\nu\cdot \nu\kappabar\Delta\neq 0\in\pi_{68,4}\TMF/\tau$, and $\nu x_3=\nu x_4=0\in\pi_{68,4}\Sph_\BP/\tau$, we must have $\eps_1=1$.

    The blue dot in bidegree $(71,5)$ indicates that, for any lift $x$ as above, $\nu^2 x$ is divisible by $4$ modulo $\tau$. Inspection of the charts shows that 
    \[\nu^2 x=4x_{71}+\tau^2(\eps_1[l_1]+\eps_2[h_2^3H_1])\]
    in $\pi_{71,5}\Sph_\BP$, for some $x_{71}$, for $[l_1]$ and $[h_2^3H_1]$ classes with $\pi([l_1])=l_1$ and $\pi([h_2^3H_1])=h_2^3H_1$ in the notation of the charts, and for $\eps_i=0,1$. We claim that $\tau^2[l_1]$ and $\tau^2[h_2^3H_1]$ go to zero in $\pi_{71,5}\TMF_\BP$, so that (up to multiplication by a unit in $\Z/8)$, the class $x_{71}$ goes to $\kappabar[\nu\Delta^2]$, completing the proof.

    Indeed, it suffices to show that $[l_1]$ and $[h_2^3H_1]$ go to zero in $\pi_{71,7}\TMF_\BP$. If not, since $\pi_{71,8}\TMF_\BP=0$, it follows that the image of $l_1$ or $h_2^3 H_1$ must be $\eta^3\kappabar\Delta^2\in\pi_{71,7}\TMF_\BP/\tau$. The class $h_2^3H_1\in\pi_{71,7}\Sph_\BP/\tau$ is divisible by $h_2$, and $\pi_{68,6}\TMF_\BP/\tau=0$, giving a contradiction in this case. In the latter case, we claim that all choices of $[l_1]$ satisfy $\tau^4\eta[l_1]=0$. Given this, it would follow that $\tau^4\eta^3[\eta\kappabar\Delta^2]=0\in\pi_{*,*}\TMF_\BP$, a contradiction, as the hidden $\eta$ extension from $(71,7)$ to $(72,10)$ and the corresponding $d_9$-differential imply that $\tau^4\eta^3[\eta\kappabar\Delta^2]\neq0$.
    
    To prove that $\tau^4\eta[l_1]=0$, note that 
\[\pi_{72,8}\Sph_\BP=\Z/2\{\tau^2[\Delta h_1h_3g^2],\tau^4[d_1g^2],[h_1l_1]\},\]
    where each $\pi([a])=a$ as before, and $[d_1g^2]$ and $[h_1l_1]$ are chosen so that $\tau^8[d_1g^2]=\tau^2[h_1l_1]=0$. It follows that
\[\tau^4\cdot\pi_{72,8}\Sph_\BP=\Z/2\{\tau^6[\Delta h_1h_3g^2]\},\]
    so if $\tau^4\eta[l_1]\neq0$, this would imply a hidden $\eta$-extension from $(71,7)$ to $(72,10)$, contradicting the $E_\infty$ chart.
\end{proof}

\begin{remark}\label{rm:etax71}
    Note that the class $x_{71}$ has the property that $\tau^{-1}(\eta x_{71})\neq0$. This follows from the hidden $\eta$ extension from $(71,5)$ to $(72,10)$ on the $E_\infty$-page of the ANSS for $\Sph$ of \cite{ANcharts}.
\end{remark}

\begin{lemma}\label{lm:71two}
    The class $x_{71}\in \pi_{71}\Sph$ admits an $\F_2$-synthetic lift $y_{71}\in \pi_{71,13}\Sph_{\F_2}$. Moreover, up to the image-of-$J$, $x_{71}$ admits a lift to the top-cell of $M(8,v_1^8)$.
\end{lemma}

\begin{proof}
    This follows from the same arguments as in \Cref{lm:47two} using fact that the $\F_2$-synthetic lift of $x_{71}$ must have filtration at least $13$ since it supports multiplication by $4$; as $47\equiv 71 \equiv 7$ modulo $8$, the image-of-$J$ arguments are analogous.
\end{proof}

\begin{proof}[Proof of \Cref{lm:71}]
    The proof of \Cref{lm:47} needs one small modification for this case. In this case, the fact that $\pi_{47,>5}\Sph_{\BP}$ consists entirely of $\tau$-torsion is replaced by the facts that $\eta\cdot\pi_{71,>5}\Sph_{\BP}$ consists entirely of $\tau$-torsion -- as shown in the proof of \Cref{lm:71one} -- and that $\eta \cdot\tau^{-1}(x_{71})\neq0$, as in \Cref{rm:etax71}. 
\end{proof}

\subsection{Lifts to the top-cell of other type $2$ complexes}\label{sssec:liftsversiontwo}
So far in our analysis, we have produced $v_2^{32}$-periodic families by lifting elements from $\Sph$ to the top-cell of $M(8,v_1^8)$, a particular type $2$-complex with a particular $v_2^{32}$-self map; see \cite[Th.1.8]{hurewicztmf}. To employ certain Toda bracket arguments, it will be critical to know that various elements produce families which are either $2$- or $4$-torsion. The arguments in this subsection are a variation on the previous subsection, and will help us to capture more of the synthetic Hurewicz image of $\TMF_\BP$.

Recall that by \cite[Th.1.1 \& Rmk.1.4]{mmmm}, the type $2$-complexes $M(2,v_1^4)$ and $M(4,v_1^4)$ have $v_2^{32}$-self maps.

\begin{lemma}\label{lm:80}
    Up to the image-of-$J$, the class $\kappabar^4$ lifts to the top-cell of $M(2,v_1^4)$.
\end{lemma}

\begin{proof}
    This follows from the same arguments as in \Cref{lm:47two} using class $g^4\in \pi_{80,16} \Sph_{\F_2}$ as the $\F_2$-synthetic lift of $\kappabar^4$.
\end{proof}

\begin{lemma}\label{lm:80choosinglift}
    Let $k\geq 0$. Then there is a class $y_k \in \pi_{80+192k, 16-r}\Sph_\BP$ which maps to $\tau^r \kappabar^4\Delta^{8k} \in \pi_{80+192k, 16-r} \TMF_\BP$ and is $2\tau^{12-r}$-torsion for some $2\leq r\leq 12$.
\end{lemma}

\begin{proof}
    The class $\tau^r \kappabar^4\Delta^{8k}$ lies in the synthetic Hurewicz image of $\TMF_\BP$ for some $r\geq 2$, as by \cite[Lm.7.17(2)]{hurewicztmf}, the class $\kappabar^4\Delta^{8k} \in \pi_{80+192k} \TMF$ lies in the classical Hurewicz image. We claim that there is a class $y_k \in \pi_{80+192k, 16-r}\Sph_\BP$ which maps to $\tau^r \kappabar^4\Delta^{8k} \in \pi_{80+192k, 16-r} \TMF_\BP$ and is $2\tau^{12-r}$-torsion for $2\leq r\leq 12$. Indeed, by \Cref{lm:80}, there is a $v_2^{32}$-periodic family of classes $x_k$ generated by $\kappabar^4$, up to the image-of-$J$, which are $2$-torsion. Let $y_k$ be a synthetic $\tau$-torsion free lift of this class $x_k$ to $\pi_{80+192k,f}\Sph_{\BP}$. As $x_k$ maps to $\kappabar^4\Delta^{8k} \in \pi_{80+192k}\TMF$, recall that the image-of-$J$ vanishes in $\TMF$, we see that $y_k$ must map to $\tau^{16-f}\kappabar^4\Delta^{8k} \in \pi_{80+192k,f}\TMF_\BP$ up to a unit, as for all $f$ this group is generated by this $\tau$-torsion free class. In particular, we see that $2\leq f\leq 16$. Moreover, as $x_k$ is $2$-torsion, then $y_k$ must be $2\tau^{f-4}$-torsion. Indeed, \emph{a priori} we see that $2y_k$ must be $\tau^i$-torsion for some $i$. This $i$ need not be greater than $f$ though, as there is no $\tau$-power torsion in $\pi_{*,<0}\Sph_\BP$. The fact that we completely understand the $0$ and $1$ line of the ANSS for $\Sph$ as well as its checkerboard pattern rule out any $\tau$-power torsion classes in filtrations $\leq 3$ in the $(80+192k)$-stem. Setting $16-r = f$ then shows that $y_k$ has the desired properties.
\end{proof}

We now have all of the tools necessary to finish our proof of \Cref{synthetichurewiczintext}.

\begin{prop}\label{lm:95}
    For all $k\geq 0$, the element $\tau^2\kappabar[\eta\Delta]^3\Delta^{8k}$ in $\pi_{95+192k,5} \TMF_\BP$, or some $\tau$-power, lies in the synthetic Hurewicz image.
\end{prop}

\begin{proof}
    Let $r$ be as in \Cref{lm:80choosinglift}. Using a Massey product computation on the $E_{13}$-page of the DSS for $\TMF$ and the synthetic Moss' theorem of \cite[Th.3.16]{smfcomputation}, one can compute
    \[\tau^r\kappabar[\eta\Delta]^3 \Delta^{8k} = \tau^{r-2}\kappabar[4\nu\Delta^3]\Delta^{8k} = \angbr{\tau^r\kappabar^4 \Delta^{8k}, 2\tau^{12-r}, \kappa} \subseteq \pi_{95+192k,5}\TMF_\BP\]
    with zero indeterminacy
    \[\kappa\pi_{81+192k,3}\TMF_\BP + \tau^r\kappabar^4\Delta^{8k} \pi_{15,r-11} \TMF_\BP = 0\]
    as $\tau^8\eta \kappa\kappabar^2 = 0$. By \Cref{lm:80choosinglift}, we have lifts of all of the elements in this Toda bracket to $\Sph_{\BP}$ and also that the associated Toda bracket in $\Sph_\BP$
    \[\angbr{y_k, 2\tau^{12-r}, \kappa} \subseteq \pi_{80+192k,5} \Sph_\BP\]
    is nonempty. In particular, by linearity of Toda bracket, any element in the above Toda bracket in the sphere is mapped to $\tau^r \kappabar^4 \Delta^{8k}$ in $\TMF_\BP$, as desired.
\end{proof}

\section{Detecting new $v_2^{32}$-periodic families in stable stems}\label{sec:periodicfamilies}
Equipped with an adequate knowledge of the synthetic Hurewicz image of $\TMF_\BP$, we can focus our attention on proving \Cref{mainthm:nonzerofamilies}. This is cut into two pieces: arguments using lower filtration of differentials\footnote{The first version of this article used a ``forcing an empty bracket'' argument here. The following filtration argument is simpler and more widely applicable.} and arguments using the Hurewicz image of $\J_0(3)$.

\subsection{Differentials with source in filtration $1$}\label{sec:todabracket}
There are a collection of classes which we know to be $\tau$-power torsion and some nonzero $\tau$-multiple lies in the synthetic Hurewicz image of $\TMF_\BP$ by \Cref{synthetichurewiczintext}. We will see that many of these classes are $\tau$-torsion free in $\Sph_\BP$ using $\J_0(3)$, but for a small collection of examples, we require the following alternative argument by realising that a differential in $\sigma(\TMF_\BP)$ cannot occur in the sphere.

\begin{theorem}\label{thm:todabracketdetection}
    Any class $x \in \pi_{s,f}\Sph_\BP$ with nonzero image in $\TMF_\BP$ given by a $\tau$-power multiple of one of the classes
    \[\nu^2\kappabar, \quad \nu[\nu\kappabar\Delta^2], \quad \kappa\kappabar^3, \quad \kappabar[\nu^2\Delta^4], \quad \kappabar[\nu\Delta^4], \quad \kappabar^3[\kappa\Delta^4], \quad \kappabar^5[\eta^2\kappabar\Delta^2]\]
    or a $\Delta^{8k}$-multiple, is $\tau$-torsion free.
\end{theorem}

The key observation used in the proof below is that we know all differentials in the ANSS for $\Sph$ with source of filtration $1$.

\begin{proof}
    First, notice that the image of $x$ in $\TMF_\BP$ is $\tau^{f-2}$-torsion and not $\tau^{f-3}$-torsion. Indeed, in $\TMF_\BP$ the class $x$ may by $\tau$-divisible, so write $y \in \pi_{s,f+n}\TMF_\BP$ such that $\tau^n y$ is the image of $x$. The class $y$ is then hit by a $d_{f+n-1}$ with source in filtration $1$ by choice of the above classes. For example, $\nu^2\kappabar$ has filtration $6$ in $\TMF_\BP$ and is hit by a $d_5$. In particular, we see that $y$ is $\tau^{f+n-2}$-torsion and not $\tau^{f+s-3}$-torsion, so the image of $x$ is $\tau^{f-2}$-torsion and not $\tau^{f-3}$-torsion, as desired. In particular, $x$ is not $\tau^{f-3}$-torsion in $\Sph_\BP$ as well. This means that the $\tau$-reduction of $x$ can only be hit by a $d_{f-1}$-differential in the ANSS for $\Sph$. However, the the $1$-line of the ANSS for $\Sph$ is well-known and cannot support a differential hitting any classes that are not $v_1$-periodic. Hence there is no room for any differentials to hit the mod $\tau$-reduction of $x$, hence $x$ is $\tau$-torsion free.
\end{proof}

\subsection{On the Hurewicz image of height $2$ image-of-$J$ spectra}\label{sec:hurewiczofjzerothree}
Combining our knowledge of the synthetic Hurewicz image of $\TMF_\BP$ together with the ``deleting differentials''-technique, we immediately obtain some elements in the classical Hurewicz image of $\J_0(3)$.

\begin{theorem}\label{hurewiczforjzerothree}
    All $v_2^{32}$-periodic families of \Cref{tab:classes} with a $\checkmark$ in the $\J_0(3)$-column have nonzero image in $\J_0(3)$.
\end{theorem}

There are many other classes in the Hurewicz image on $\J_0(3)$ in the image of the boundary map $\pi_{\ast+1} \TMF_0(3) \to \pi_\ast \J_0(3)$. These include all of the divided $\al$-family and any $v_2$-periodic families will be explored in future work together with the Hurewicz image of $\Q(3)$.

\begin{proof}
    Consider \Cref{cor:synthetichurewiczimpractice} with
    \[R=\Sph_\BP, \qquad A=\TMF_\BP, \qquad B=\TMF_0(3), \qquad \varphi = q, \qquad \psi = p\]
    and $x=\nu\kappabar$. In this case, we have $d_5$-differentials $d_5(\pm(4k+1)\Delta) = \nu\kappabar$ for all integers $k$ by \cite[Pr.6.19]{smfcomputation}. These are all possible differentials witnessing $\nu\kappabar$ as a $d_5$-boundary, up to a sign. Lifting these elements to $\TMF_\BP/\tau^4$ as in \cite[Not.6.18]{smfcomputation}, we see that $(q-p)(\Delta)$ is nonzero courtesy of \Cref{pr:injectivityofpandq}. In particular, $\nu\kappabar$ is not a $d_5$-boundary in $\sigma(\J_0(3)_\BP)$ by \Cref{cor:synthetichurewiczimpractice}. The class $\nu\kappabar$ lives in $\pi_{23,5}\J_0(3)_\BP$ and $\pi_{n,f}\J_0(3)_\BP/\tau$ vanishes for $f\leq -1$ by \Cref{lm:nothingbelowzero}, so $\tau$-inverts to a nonzero class of the same name in the Hurewicz image of $\J_0(3)$ by \Cref{cor:synthetichurewiczimpractice}. Moreover, the $2$-extension $4\nu = \tau^2 \eta^3$, see \cite[Lm.6.16]{smfcomputation}, shows that $\eta^3\kappabar$ lies in the Hurewicz image of $\J_0(3)$ too, and that $\nu\kappabar$ has order $8$. This argument works equally as well for all elements in the $v_2^{32}$-periodic family generated by $\nu\kappabar$.

    Every other family in \Cref{tab:classes} with a $\checkmark$ in the $\J_0(3)$-column lies in the Hurewicz image of $\J_0(3)$ by the same argument. To summarise these, we have the following table indicating the necessary differentials; references to all of these differentials can be found on their appropriate pages in \cite[\textsection6]{smfcomputation}. We only need check one generator in each degree, as other generators in higher or lower filtration are connected by extensions, such as $\tau^2 \eta^3 = 4\nu^2$ in $\pi_{3,1}\TMF_\BP$ or $\nu[\eta\kappabar\Delta] = \tau^4 \eps\kappabar^2$ in $\pi_{48,6}\TMF_\BP$, all of which are displayed in \cite{bauer} or \cite{konter}, and proven in either \cite{bauer} or \cite{smfcomputation}.

    \begin{table}[h]
\centering
{
\begin{tabular}{|c|c|c|c|c|c|}
\hline
Degree  &   AN-filt.&Group  &   Target                                &   $E_2$-rep.        &   Differential                                            \\  \hline
$47$    &   $5$     &$\Z/4$ &   ${[2\nu\kappabar\Delta]}$             &   $2h_2g\Delta$     &   $d_5((2k+1)\Delta^2) = 2\nu\kappabar\Delta$             \\  
$48$    &   $10$    &$\Z/2$ &   ${\eps\kappabar^2}$                   &   $cg^2$            &   $d_9([\eta\Delta^2]) = \eps\kappabar^2$                 \\  

$71$    &   $5$     &$\Z/8$ &   ${[\nu\kappabar\Delta^2]}$            &   $h_2g\Delta^2$    &   $d_5((4k+1)\Delta^3) =\pm\nu\kappabar\Delta^2$          \\  
$72$    &   $10$    &$\Z/2$ &   $\kappabar^2[\eps\Delta]$             &   $cg\Delta^2$      &   $d_9([\eta\Delta^3]) = \kappabar^2[\eps\Delta]$         \\  
$73$    &   $11$    &$\Z/2$ &   $\eta{\kappabar[\eps\kappabar\Delta]}$&   $h_1cg\Delta^2$   &   $d_9([\eta^2\Delta^3]) = \eta\kappabar^2[\eps\Delta]$   \\  

$95$    &   $7$     &$\Z/2$ &  ${[\eta^3\kappabar\Delta^3]}$          &   $h_1^3g\Delta^3$  &   $d_7((2k+1)\Delta^4) = \kappabar[\eta\Delta]^3$         \\  

$119$   &   $5$     &$\Z/4$ &   $\nu{[2\kappabar\Delta^4]}$           &   $2h_2g\Delta^4$   &   $d_5((4k+2)\Delta^5) = \nu[2\kappabar\Delta^4]$         \\  
$120$   &   $24$    &$\Z/2$ &   ${\kappabar^6}$                       &   $g^6$             &   $d_{23}([\eta\Delta^5]) = \kappabar^6$                  \\  

$143$   &   $7$     &$\Z/2$ &   $\eta{[\eta^2\kappabar\Delta^5]}$     &   $h_1^3g\Delta^5$  &   $d_5(\Delta^6) = \eta\kappabar\eta\Delta^6$             \\  
$144$   &   $10$    &$\Z/2$ &   ${\kappabar^2[\eps\Delta^4]}$         &   $cg^2\Delta^4$    &   $d_9([\eta\Delta^6]) = \kappabar^2[\eps\Delta^4]$       \\  
$145$   &   $23$    &$\Z/2$ &   ${\kappabar^5[\eta\kappabar\Delta]}$  &   $h_1g^6\Delta$    &   $d_{23}([\eta^2\Delta^6]) = \kappabar^6[\eta\Delta]$    \\  

$167$   &   $7$     &$\Z/2$ &   $\kappabar{[\eta^3\Delta^6]}$         &   $h_1^3g\Delta^6$  &   $d_7([(8k+4)\Delta^7]) = \kappabar[\eta^3\Delta^6]$     \\  
$168$   &   $10$    &$\Z/2$ &   ${\kappabar^2[\eps\Delta^5]}$         &   $cg^2\Delta^5$    &   $d_9([\eta\Delta^7]) = \kappabar^2[\eps\Delta^5]$       \\  
$169$   &   $11$    &$\Z/2$ &   $\eta{\kappabar^2[\eps\Delta^5]}$     &   $h_1cg^2\Delta^5$ &   $d_9(\eta[\eta\Delta^7])=\eta\kappabar^2[\eps\Delta^5]$ \\  \hline
\end{tabular}
}
\caption{
{Table of differentials in $\sigma(\TMF_\BP)$ needed to apply \Cref{cor:synthetichurewiczimpractice}.
}}
    \label{tab:differentials}
\end{table}

The target of each of these differentials vanishes under $q-p$ as they all lie in the synthetic Hurewicz image of $\TMF_\BP$ or $\tau^2$ is an isomorphism into this synthetic Hurewicz image. Moreover, by \Cref{pr:injectivityofpandq} and inspection of $\sigma(\TMF_0(3)_\BP)$, so the ANSS for $\TMF_0(3)$, we see that the sources of all of these differentials have nonzero image under $q-p$. Note that the degree of each $d_r$-differential is either equal to the AN-filtration of its target, or at most smaller by 2. In particular, not only do these families lift to elements in $\sigma(\J_0(3)_\BP)$ which are not $d_r$-boundaries, but using \Cref{lm:checkboard,lm:nothingbelowzero,lm:filtrationzero}, they are all permanent cycles, hence represent classes in the classical Hurewicz image of $\J_0(3)$. This finishes the proof.
\end{proof}

\begin{remark}
    The families in degrees congruent to $23, 47, 71, 95, 119, 143$, and $167$ also lie in the Hurewicz image of the fibre of $\psi^3-1\colon \TMF\to \TMF$. The key fact is that $\psi^3-1$ acts on $\pi_{\ast,\ast}\TMF_\BP$ injectively in AN-filtration $0$ and zero elsewhere. In fact, this holds for the fibre of $\psi^N-1$ for any odd integer $N$. We leave the details to the reader.
\end{remark}

\subsection{Collecting all one hundred and twenty-five $v_2^{32}$-periodic families}\label{ssec:countingmethodology}
Finally, everything is in place to prove \Cref{mainthm:nonzerofamilies}, which we restate now for the readers convenience.

\begin{theorem}[{\Cref{mainthm:nonzerofamilies}}]
    The classes in \Cref{tab:classes} together represent $125$ nonzero $v_2^{32}$-periodic families in $\pi_\ast\Sph$ of order indicated by the Group-column and all with trivial image in $\pi_\ast\TMF$.
\end{theorem}

\begin{proof}
    By \Cref{synthetichurewiczintext}, all of the generators of \Cref{tab:classes} admit unique $\BP$-synthetic lifts to classes in $\pi_{\ast,\ast}\TMF_\BP$, and some nonzero $\tau$-power multiple of these classes lies in the synthetic Hurewicz image of $\TMF_\BP$. In particular, we have a collection of $v_2^{32}$-periodic families in $\pi_{\ast,\ast}\Sph_\BP$ which map to $\tau$-power torsion elements in $\pi_{\ast,\ast}\TMF_\BP$. By \Cref{thm:todabracketdetection} and \Cref{hurewiczforjzerothree}, we know that all of these classes in $\Sph_\BP$ are in fact $\tau$-torsion free and hence represent nontrivial elements in $\pi_\ast\Sph$.

    All that remains is to count these families. To this end, we will decompose each generator of \Cref{tab:classes} as an element in $\TMF_\BP$ into as many $v_2^{32}$-periodic factors as we can. For those classes coming from \Cref{thm:todabracketdetection}, this works as per the statement of \Cref{thm:todabracketdetection}. Otherwise, we note that for those generators of \Cref{tab:classes} with a $\checkmark$ in the $\J_0(3)$-column all lift to $\tau$-torsion free elements in $\J_0(3)_\BP$, this implies that all of these $v_2^{32}$-periodic factors as nonzero in $\pi_\ast \J_0(3)$.\footnote{Of course, different $v_2^{32}$-periodic families mapped to the same generator in $\TMF_\BP$ may be related, see \Cref{cortoinfiniteagreement}, but they still may define distinct families; we refer the reader to our discussion in \Cref{ssec:namingperiodicfamilies} for some more discussion.}
    
\begin{itemize}
    \item[23$a$)]    The periodicity of $\nu\kappabar$ comes from $\kappabar$ by \cite[Lm.7.17(2)]{hurewicztmf}. This yields $1$ nonzero $v_2^{32}$-periodic family.
    \item[23$b$)]    In $\pi_{23,5} \TMF_\BP$ we have the equalities
    \[4\nu\kappabar = \tau^2\eta^3\kappabar = \eta \eps \kappa = \nu^3 \kappa,\]
    which follow from the relations $4\nu = \tau^2 \eta^3$ in $\pi_{3,1} \TMF_\BP$, $\eps\kappa = \tau^2 \eta^2\kappabar$ in $\pi_{22,4} \TMF_\BP$, and $\nu^3 = \eta\eps\in \pi_{9,3} \TMF_\BP$, the first of which holds in $\Sph_\BP$, then second holds from the Toda bracket manipulation
    \[\tau^2 \eta^2 \kappabar = \eta \angbr{\nu, 2\nu, \kappa} = \angbr{\eta, 2\nu, \nu}\kappa = \eps \kappa\]
    using \cite[Pr.5.11]{smfcomputation}, and the third from an $E_2$-page computation of $\sigma(\TMF_\BP)$. These equalities also hold in $\pi_{23}\Sph$ after $\tau$-inversion; see \cite[Th.11.61]{brunerrognes}. The classes $\nu^2$, $\eps$, and $\kappa$ are all $v_2^{32}$-periodic in $\Sph$ by \cite[Lms.7.17(1) \& 7.27]{hurewicztmf}; note that $\nu^2$ and $\eps$ are not shown to lift to the top-cell of $M(8,v_1^8)$, but rather that $\nu^2\Delta^8$ and $\eps\Delta^8$ do.\footnote{We do not count the periodicity from $\kappabar$ again, as this family is simply the $2$-torsion in the family generated by $\nu\kappabar$. This remark also holds in the degrees $47b$, $71b$, $74b$, $119b$, $170b$, and $170c$.} This yields $4$ nonzero $v_2^{32}$-periodic families.
    \item[26)]     Using \Cref{thm:todabracketdetection}, the periodicity of $\nu^2\kappabar$ comes from both $\nu^2$ and $\kappabar$, which we have seen are $v_2^{32}$-periodic. This yields $2$ nonzero $v_2^{32}$-periodic families.
    \item[47$a$)]    The periodicity of $[2\nu\kappabar\Delta]$ comes from \Cref{lm:47}. This yields $1$ nonzero $v_2^{32}$-periodic family.
    \item[47$b$)]    In $\pi_{47,5} \TMF_\BP$ we have the equalities
    \[2[2\nu\kappabar\Delta] = \tau^2\eta^2 [\eta \kappabar\Delta] = \eps [\eta \kappa\Delta] = \eta \kappa [\eps\Delta]\]
    which follow from the equalities used in degree 23b in $\TMF_\BP/\tau^4$. For degree reasons these relations then persist to $\TMF_\BP$. We have seen that $\eps$ and $\kappa$ are $v_2^{32}$-periodic, and so are $[\eps\Delta]$, $[\eta\kappa\Delta]$, and $[\eta\kappabar\Delta]$ by \cite[Lm.7.18(1)-(3)]{hurewicztmf}. This yields $5$ nonzero $v_2^{32}$-periodic families.
    \item[48)]     There is the equality $\eps\kappabar = \kappa^2$ in $\pi_{28}\Sph$, see \cite[Th.11.61(28)]{brunerrognes}, and since $\pi_{28,6} \TMF_\BP$ is $\tau$-torsion free, the relation $\kappa^2 = \tau^2 \eps \kappabar$ in $\TMF_\BP$ is forced. This is the third equality in $\pi_{48,6} \TMF_\BP$
    \[\eta [2\nu\kappabar\Delta] = \nu [\eta\kappabar\Delta] = \tau^2 \kappa^2 \kappabar = \tau^4 \eps \kappabar^2,\]
    the second is the $\nu$-extension from \cite[Lm.6.40(28,2)]{smfcomputation} and the first is similar. We have seen that $\eps$, $\kappa$, $\kappabar$, $[\eta\kappabar\Delta]$, and $[2\nu\kappabar\Delta]$ are all $v_2^{32}$-periodic. 
    This seems to yield $6$ nonzero $v_2^{32}$-periodic families, but we have over-counted by $1$. Indeed, the classes $\kappa^2\kappabar$ and $\eps\kappabar^2$ generate two periodic families from the $v_2^{32}$-periodicity of $\kappabar$, but these families are the same, as in $\pi_{28}\Sph$ we have $\kappa^2 = \eps\kappabar$ by \cite[Th.11.61(28)]{brunerrognes}. This then yields a total of $5$ nonzero $v_2^{32}$-periodic families.
    \item[71$a$)]    The periodicity of $[\nu\kappabar\Delta^2]$ comes from \Cref{lm:71}. This yields $1$ nonzero $v_2^{32}$-periodic family.
    \item[71$b$)]    As in degrees $23b$ and $47b$, there are equalities in $\pi_{71,5} \TMF_\BP$
    \[4[\nu\kappabar\Delta^2] = \tau^2 \eta [\eta^2 \kappabar\Delta^2] = [\eps\Delta][\eta\kappa\Delta] = \nu \kappa [\nu^2\Delta^2] = \nu^2[\nu\kappa \Delta^2].\]
    We have seen that $\nu^2$, $\kappa$, $[\eps\Delta]$, and $[\eta\kappa\Delta]$ are $v_2^{32}$-periodic, and so are $[\nu^2\Delta^2]$, $[\nu\kappa\Delta^2]$, and $[\eta^2\kappabar\Delta^2]$ by \cite[Lm.7.21(1)-(3)]{hurewicztmf}. This yields $7$ nonzero $v_2^{32}$-periodic families.
    \item[72)]     The extension $\eta[\nu\Delta^2] = \tau^4 \kappabar[\eps\Delta]$ of \cite[Lm.6.40(28,2)]{smfcomputation} mentioned above also implies that $\eta [\nu\kappabar\Delta^2] = \tau^4 \kappabar^2[\eps\Delta]$ in $\pi_{72,6} \TMF_\BP$. We have seen that $\kappabar$, $[\eps\Delta]$, and $[\nu\kappabar\Delta^2]$ are $v_2^{32}$-periodic. This yields $3$ nonzero $v_2^{32}$-periodic families.
    \item[73)]     In $\pi_{73,7} \TMF_\BP$ we have equalities
    \[\eta^2 [\nu\kappabar\Delta^2] = \nu [\eta^2\kappabar\Delta^2] = \tau^4 \eta \kappabar^2 [\eps\Delta] = \tau^4 \eps \kappabar [\eta\kappabar\Delta] = \tau^4 \kappa\kappabar [\eta \kappa\Delta],\]
    using the algebra of the $E_2$-page of the ANSS for $\TMF$ together with similar extensions used in degree $72$.
    We already know that $\eps$, $\kappa$, $\kappabar$, $[\eps\Delta]$, $[\eta\kappa\Delta]$, $[\eta\kappabar\Delta]$, $[\eta^2 \kappabar\Delta^2]$, and $[\nu \kappabar\Delta^2]$ are all $v_2^{32}$-periodic. This yields $10$ nonzero $v_2^{32}$-periodic families.
    \item[74$a$)]    Using \Cref{thm:todabracketdetection}, the periodicity of $\nu[\nu\kappabar\Delta^2] = \kappabar[\nu^2\Delta^2]$, an equality which clearly holds on $E_2$-pages and lifts to $\TMF_\BP$ as there are no classes in higher filtration, comes from $\kappabar$, $[\nu^2\Delta^2]$, and $[\nu\kappabar\Delta^2]$, which we have seen are $v_2^{32}$-periodic. This yields $3$ nonzero $v_2^{32}$-periodic families.
    \item[74$b$)]    Using \Cref{thm:todabracketdetection}, the periodicity of
    \[2\nu[\nu\kappabar\Delta^2] = \tau^8 \kappa\kappabar^3,\]
    an extension which follows from \cite[Lm.6.53(54,2)]{smfcomputation}, comes from both $\kappa$ and $\kappabar$, which we have seen are $v_2^{32}$-periodic. This yields $2$ nonzero $v_2^{32}$-periodic families.
    \item[95)]     The periodicity of $[4\nu\kappabar\Delta^3]$ comes from \Cref{lm:95}.\footnote{The infinite $v_2^{32}$-periodic family in this congruence class of degrees is not necessarily generated by a single class, but rather by a collection of $v_2^{32}$-periodic nonempty Toda brackets which do not contain zero; see the proof of \Cref{lm:95}.} This yields $1$ nonzero $v_2^{32}$-periodic family.
    \item[119$a$)]   The periodicity of $\nu[2\kappabar\Delta^4]$ comes from that of $[2\kappabar\Delta^4]$ shown in \cite[Lm.7.22(2)]{hurewicztmf}. This yields $1$ nonzero $v_2^{32}$-periodic family.
    \item[119$b$)]   In $\pi_{119,5} \TMF_\BP$ we have the equalities
    \[2\nu [2\kappabar\Delta^4] = \tau^2 \eta^2 [\eta\kappabar\Delta^4] = \eta\kappa[\eps\Delta^4] = \eta \eps [\kappa\Delta^4] = \nu^3 [\kappa\Delta^4] = \nu\kappa[\nu^2\Delta^4] = [\nu^2\Delta^2][\nu\kappa\Delta^2],\]
    as these hold in $\TMF_\BP/\tau^4$ and lift to $\TMF_\BP$, similar to the arguments in degree $23b, 47b$, and $71b$. We know that $\nu^2$, $\eps$, $\kappa$, $[\nu^2\Delta^2]$, and $[\nu\kappa\Delta^2]$ are $v_2^{32}$-periodic, and so are $[\nu^2\Delta^4]$, $[\eps\Delta^4]$, $[\kappa\Delta^4]$, and $[\eta\kappabar\Delta^4]$ by \cite[Lms.7.22(1) \& 7.23]{hurewicztmf}. This yields $11$ nonzero $v_2^{32}$-periodic families.
    \item[120)]    There is an equality $\tau^{18} \kappabar^6 = \nu [\eta\kappabar\Delta^4]$ of \cite[Rmk.8.3]{smfcomputation} in $\pi_{120,6} \TMF_\BP$. We already know that $\kappabar$ and $[\eta\kappabar\Delta^4]$ are $v_2^{32}$-periodic. This yields $2$ nonzero $v_2^{32}$-periodic families.
    \item[122)]    Using \Cref{thm:todabracketdetection}, the periodicity of $\kappabar[\nu^2\Delta^4]$ comes from both $\kappabar$ and $[\nu^2\Delta^4]$, which we have seen are $v_2^{32}$-periodic. This yields $2$ nonzero $v_2^{32}$-periodic families.
    \item[143)]    In $\pi_{143,5} \TMF_\BP$ we have the equalities
    \[\eps[\eta\kappa\Delta^5] = \eta \kappa[\eps\Delta^5] = \eta [\eps\Delta][\kappa\Delta^4] = [\eta\kappa\Delta][\eps\Delta^4],\]
    as this is true on $E_2$-pages and lifts due to lack of elements in higher filtration. We know that $\eps$, $\kappa$, $[\eps\Delta]$, $[\eta\kappa\Delta]$, $[\eps\Delta^4]$, and $[\kappa\Delta^4]$ are $v_2^{32}$-periodic, and the periodicity of $[\eps\Delta^5]$ and $[\eta\kappa\Delta^5]$ is shown in \cite[Lms.7.24 \& 7.25]{hurewicztmf}. This yields $8$ nonzero $v_2^{32}$-periodic families.
    \item[144)]    In $\pi_{144,8} \TMF_\BP$ we have an equality
    \[\tau^2 \kappabar^2 [\eps\Delta^4] = \kappa\kappabar [\kappa\Delta^4].\]
    We have seen that $\kappa$, $\kappabar$, $[\eps\Delta^4]$, and $[\kappa\Delta^4]$ are $v_2^{32}$-periodic. This yields $5$ nonzero $v_2^{32}$-periodic families.
    \item[145)]    In $\pi_{145,5} \TMF_\BP$ we have the equalities
    \[\nu\kappa[\eps\Delta^5] = \tau^6 \eta\kappabar^2[\eps\Delta^4] = \tau^6 \eps \kappabar [\eta \kappabar \Delta^4] = \tau^4 \eta \kappa \kappabar[\kappa\Delta^4] = \nu [\eps\Delta][\kappa\Delta^4] = \tau^{20} \kappabar^5[\eta\kappabar\Delta],\]
    forced by $E_2$-page computations, known extensions, and the last equality is an extension drawn in \cite{bauer}, for example, which one obtains synthetically using the total differential technique of \cite[Rmk.8.3]{smfcomputation}. We know that $\eps$, $\kappa$, $\kappabar$, $[\eps\Delta]$, $[\eta\kappabar\Delta]$, $[\eps\Delta^4]$, $[\kappa\Delta^4]$, $[\eta\kappabar\Delta^4]$, and $[\eps\Delta^5]$ are $v_2^{32}$-periodic. By \cite[Th.11.61(35)]{brunerrognes}, we have $\nu[\eps\Delta] = \eta\kappa\kappabar$ in the sphere, so the two families generated by $[\kappa\Delta^4]$ above are the same. This seems to yield $14$ nonzero $v_2^{32}$-periodic families, but we have over-counted by $1$. Indeed, using the classes $\eta\kappa\kappabar[\kappa\Delta^4]$ and $\nu[\eps\Delta][\kappa\Delta^4]$ to generate families from the $v_2^{32}$-periodicity of $[\kappa\Delta^4]$ give the same families, as $\nu[\eps\Delta] = \eta\kappa\kappabar$ in $\pi_{35} \Sph$ by \cite[Th.11.61(35)]{brunerrognes}. This then yields $13$ total nonzero $v_2^{32}$-periodic families.
    \item[167)]    In $\pi_{167,5} \TMF_\BP$ we have the equalities
    \[\nu^2 [\nu\kappa\Delta^6] = \nu \kappa [\nu^2 \Delta^6] = \nu [\nu^2\Delta^2][\kappa \Delta^4] = [\nu\kappa\Delta^2][\nu^2\Delta^4]= [\eps\Delta][\eta\kappa\Delta^5] = [\eta\kappa\Delta][\eps\Delta^5]\]
    all from lifting $E_2$-page computations. We know that $\nu^2$, $\kappa$, $[\eps\Delta]$, $[\eta\kappa\Delta]$, $[\nu^2\Delta^2]$, $[\nu\kappa\Delta^2]$, $[\nu^2\Delta^4]$, $[\kappa\Delta^4]$, $[\eps\Delta^5]$, and $[\eta\kappa\Delta^5]$ are $v_2^{32}$-periodic, and $[\nu^2\Delta^6]$ and $[\nu\kappa\Delta^6]$ are $v_2^{32}$-periodic by \cite[Lm.7.26(1)-(2)]{hurewicztmf}. This yields $12$ nonzero $v_2^{32}$-periodic families.
    \item[168)]    The periodicity of $\kappabar^2[\eps\Delta^5]$ comes from both $\kappabar$ and $[\eps\Delta^5]$. This yields $2$ nonzero $v_2^{32}$-periodic families.
    \item[169)]    In $\pi_{169,9} \TMF_\BP$ we have the equalities
    \[\kappa\kappabar [\eta \kappa\Delta^5]= \kappabar [\eta\kappa\Delta][\kappa\Delta^4] =\tau^2 \eta\kappabar^2 [\eps\Delta^5] = \tau^2 \kappabar[\eps\Delta][\eta \kappabar\Delta^4] = \tau^2 \kappabar[\eta\kappabar\Delta] [\eps\Delta^4]\]
    from $E_2$-page computations and the relation $\kappa^2 = \tau^2\eps\kappabar$. We know that $\kappa$, $\kappabar$, $[\eps\Delta]$, $[\eta\kappa\Delta]$, $[\eta\kappabar\Delta]$, $[\eps\Delta^4]$, $[\kappa\Delta^4]$, $[\eta\kappabar\Delta^4]$, $[\eps\Delta^5]$, and $[\eta\kappa\Delta^5]$ are all $v_2^{32}$-periodic. This yields $14$ nonzero $v_2^{32}$-periodic families.
    \item[170$a$)]   Using \Cref{thm:todabracketdetection}, the periodicity of $\kappabar[\nu^2\Delta^6]$ comes from $\kappabar$ and $[\nu^2\Delta^6]$, which we have seen are $v_2^{32}$-periodic. This yields $2$ nonzero $v_2^{32}$-periodic families.
    \item[170$b$)]   Using \Cref{thm:todabracketdetection}, the periodicity of
    \[2\kappabar[\nu^2\Delta^6] = \tau^8 \kappabar^3 [\kappa\Delta^4],\]
    where this equality is given in \cite[Lm.6.53(150,2)]{smfcomputation}, comes from $\kappabar$ and $[\kappa\Delta^4]$, which we have seen are $v_2^{32}$-periodic. This yields $2$ nonzero $v_2^{32}$-periodic families.
    \item[170$c$)]   Using \Cref{thm:todabracketdetection}, the periodicity of
    \[4\kappabar[\nu^2\Delta^6] = \nu^3[\nu\kappa\Delta^6] = \tau^6 \eta^2 \kappabar^2 [\eps\Delta^5] = \tau^{20} \kappabar^5[\eta^2\kappabar\Delta^2],\]
    where the last equality is the extension proven in \cite[Cor.C]{smfcomputation}, comes from $\nu^2$, $\kappabar$, $[\eta^2\kappabar\Delta^2]$, $[\eps\Delta^6]$, and $[\nu\kappa\Delta^6]$, which we have seen are $v_2^{32}$-periodic. This yields $6$ nonzero $v_2^{32}$-periodic families.
\end{itemize}
\end{proof}

\begin{remark}\label{rmk:televsk2}
    All of the families produced by \Cref{mainthm:nonzerofamilies} have nonzero image in $\Sph_{T(2)}$ by definition. In fact, we claim that actually all elements in these families have nontrivial image in $\Sph_{K(2)}$, the coarser localisation of the sphere spectrum, now with respect to the second Morava $K$-theory spectrum. Ravenel's telescope conjecture asked if the canonical map $\Sph_{T(h)} \to \Sph_{K(h)}$ is an equivalence for all $h\geq 1$ at all primes; a statement decisively proven in the negative in \cite{telescope} for all $h\geq 2$ and at all primes. Nevertheless, our techniques are based on the ANSS for $\TMF$, which in positive filtration agrees with the $G_{48}$-HFPSS for $E_2$ converging to the $K(2)$-localisation of $\TMF$; see \cite{g48}. It follows that all of the families of \Cref{mainthm:nonzerofamilies} have nontrivial image in the $K(2)$-localisation of the sphere. In fact, as $\TMF$, and hence also $\J_0(3)$, are $\MU$-nilpotent, and $\MU$-nilpotent spectra satisfy the telescope conjecture, there is no chance that these techniques could be used to uncover nonzero families in the $T(2)$-local sphere which vanish $K(2)$-locally.

\end{remark}

\subsection{Further questions}\label{ssec:questions}
The deleting differentials and forcing an empty Toda bracket techniques were pushed as far as they could be here using $\J_0(3)$ and given the state of the synthetic Hurewicz image of $\TMF_\BP$ in \Cref{synthetichurewiczintext}. If any of these variables are changed or improved upon, this could lead to further new $v_2^{32}$-periodic families. Here we pose some of these questions, as well as some suggested techniques.

\subsubsection{Synthetic Hurewicz image of topological modular forms}
The arguments of \Cref{lm:47,lm:71,lm:95} are clearly restricted by how much of the ANSS for $\Sph$ we know, i.e., how large the charts in \cite{ANcharts}. For a more systematic approach, one might consider a $\BP$-synthetic version of the arguments made by Behrens--Mahowald--Quigley in \cite{hurewicztmf} using a modified ASS internal to $\BP$-synthetic spectra to fully compute the synthetic Hurewicz image of $\TMF_\BP$.

\begin{question}\label{question:theselementsinsynthetitmf}
    Are the generators
    \[\kappabar[\nu\Delta^4]\in \pi_{119,5} \TMF_\BP, \qquad
    \kappabar[2\nu\Delta^5]\in \pi_{143,5} \TMF_\BP, \qquad
    \kappabar[\nu\Delta^6]\in \pi_{167,5} \TMF_\BP,\]
    or any nonzero $\tau$-power multiples, in the synthetic Hurewicz image? Do there exists lifts of these classes to $\Sph_\BP$ which are $\tau$-torsion free and are $v_2^{32}$-periodic? Could this be approached using longer synthetic Toda brackets such as
    \[\angbr{\kappabar^5, \tau^4\nu, \tau^4 2\nu, \tau^4 3\nu, \tau^4 4\nu, 5\nu} \subseteq \pi_{119,5} \TMF_\BP?\]
\end{question}

\begin{question}
    Do the classes $\nu^2, [\nu^2\Delta^2], [\nu^2\Delta^4]$, and $[\nu^2\Delta^6]$ form a $v_2^8$-periodic family as suggested by \cite[\textsection13]{tmfbook}? What about for the related families starting with $\nu^3$ and $\nu\kappa$? Is there a simple type $2$ finite complex smaller than those $\mathcal{Z}$ of Bhattacharya--Egger \cite{v21selfmaps} which witnesses this?
\end{question}

\begin{question}
    Do the classes $\eps, [\eps\Delta^4]$ form a $v_2^{16}$-periodic family? What about related pairs such as $\kappa,[\kappa\Delta^4]$ and $[\eps\Delta], [\eps\Delta^5]$ and various multiples?
\end{question}

\begin{remark}
    If one knew more about the synthetic Hurewicz image of $\TMF_\BP$, then one could hope to improve further upon \Cref{mainthm:nonzerofamilies} utilising \Cref{cor:synthetichurewiczimpractice}. For example, if all of the classes of \Cref{question:theselementsinsynthetitmf} above lie in the synthetic Hurewicz image of $\TMF_\BP$, then their $\tau$-inversions would lift to elements in the Hurewicz image of $\J_0(3)$ from the same arguments made in the proof of \Cref{hurewiczforjzerothree} above.
\end{remark}

\subsubsection{Using other detection spectra}
One can try to push the computations of \Cref{hurewiczforjzerothree} further using other versions of $\J_0(N)$ or Behrens' $\Q(N)$-spectra of \cite{ktwospheremark}.

\begin{question}\label{questionnumberthesecond}
    Can one apply \Cref{cor:synthetichurewiczimpractice} to $\J_0(5)$ at $p=2$ or $\J_0(7)$ at $p=3$, using the computations of \cite{markkyle,meieroroztmfo7}, to determine any more nonzero $v_2^{32}$- or $v_2^9$-periodic families? What about fibres of Adams operations on other periodic spectra, such as fixed points of Lubin--Tate theories (\cite{ktwolocalsphere}) or Barsotti--Tate theories (see \cite[\textsection6.1]{luriestheorem}), or the Adams operations on connective topological modular forms from \cite{realspectra,adamsontmf}?
\end{question}

\begin{question}\label{questionrelationtoQN}
    Does an iterative version of \Cref{cor:synthetichurewiczimpractice} apply to synthetic versions of Behrens' $\Q(N)$-spectra, given by computing a Bousfield--Kan spectral sequence internal to synthetic spectra, to yield even more families, such as a potential $2$-torsion family generated by $\kappabar[\nu^2\Delta^4]$ in degree $122$ or the $8$-torsion family generated by $\kappabar[\nu^2\Delta^6]$ in degree $170$?
\end{question}

In low degrees, it seems that $\Q(3)$ also captures the $2$-torsion family generated by $\nu^2\kappabar$ in degree $26$ and the $4$-torsion family generated by $\kappabar[\nu^2\Delta^2]$ in degree $74$. Of course, we have already seen that these classes generated nonzero $v_2^{32}$-periodic families by \Cref{thm:todabracketdetection}.

\begin{remark}
    There is a somewhat orthogonal method to detecting periodic families in $\J_0(3)$ and $\Q(3)$ which we have avoided in this article for simplicity. Arguing as we did at height $1$ in \cite{syntheticj} or height $2$ at the prime $3$ in \cite{heighttwojat3}, one can use $\Q(3)_\BP = \lim \nu\Q(3)^\bullet$, where $\Q(3)^\bullet$ is the diagram defining $\Q(3)$, to detect \emph{divided $\beta$-family} elements. This method, for example, can be used to recover many of the computations of Behrens--Ormsby \cite{markkyle}. For example, that $\Q(3)$ detects $\sigma^2 = \be_{4/4}$. We will return to this in future work.
\end{remark}



\scriptsize
\bibliography{references} 
\bibliographystyle{alpha}

\end{document}